\documentclass[a4paper,11pt,reqno]{amsart}
\usepackage{mhequ}
\usepackage{enumerate}
\usepackage{amssymb, amsmath}
\usepackage{mathrsfs,dsfont}
\usepackage{amscd}
\usepackage[active]{srcltx}
\usepackage{verbatim}
\usepackage[colorlinks,linkcolor={blue},citecolor={blue},urlcolor={black}]{hyperref}
\usepackage[]{changebar}

\let\emptyset \undefined
\newsymbol\emptyset    203F

\theoremstyle{plain}
\newtheorem{theorem}{Theorem}[section]
\newtheorem{corollary}[theorem]{Corollary}
\newtheorem{lemma}[theorem]{Lemma}
\newtheorem{proposition}[theorem]{Proposition}

\newtheorem{definition}[theorem]{Definition}
\newtheorem{assumption}[theorem]{Assumption}

\newtheorem*{definition*}{Definition}

\theoremstyle{remark}
\newtheorem{remark}[theorem]{Remark}
\newtheorem{example}[theorem]{Example}
\newtheorem*{remark*}{Remark}
\newtheorem*{example*}{Example}
\newtheorem*{notation*}{Notation}

\numberwithin{equation}{section}

\def\R{{\mathbb R}}

\newcommand{\F}{{\mathcal F}}

\newcommand{\eps}{\varepsilon}
\renewcommand{\phi}{\varphi}

\newcommand{\dd}{\; \mathrm{d}}

\DeclareMathOperator{\supp}{supp}

\DeclareMathOperator{\Ran}{Ran}

\DeclareMathOperator{\Ker}{Ker}

\DeclareMathOperator{\sgn}{sgn}

\newcommand{\ip}[1]{\langle {#1}\rangle}
\DeclareMathOperator{\arctanh}{arctanh}
\DeclareMathOperator{\diag}{diag}
\DeclareMathOperator{\grad}{grad}

\DeclareMathOperator{\Ric}{Ric}

\newcommand{\beq}{\begin{equation}}
\newcommand{\eeq}{\end{equation}}
\newcommand{\bal}{\begin{aligned}}
\newcommand{\eal}{\end{aligned}}
\newcommand{\ben}{\begin{enumerate}}
\newcommand{\beni} {\begin{enumerate}[(i)]}
\newcommand{\een}{\end{enumerate}}
\newcommand{\bit}{\begin{itemize}}
\newcommand{\eit}{\end{itemize}}
\newcommand{\beqw}{\begin{equation*}}
\newcommand{\eeqw}{\end{equation*}}
\newcommand{\bex}{\begin{example}}
\newcommand{\eex}{\end{example}}
\newcommand{\bre}{\begin{example}}
\newcommand{\ere}{\end{example}}
\newcommand{\bma}{\begin{bmatrix}}
\newcommand{\ema}{\end{bmatrix}}

\newcommand{\Dom}{{\mathsf D}}

\newcommand{\one}{{{\bf 1}}}

\renewcommand{\hat}{\widehat}

\newcommand{\ddt}{\frac{\mathrm{d}}{\mathrm{d}t}}

\newcommand{\cQ}{\mathcal{Q}}

\newcommand{\cM}{\mathcal{M}}
\newcommand{\cH}{\mathcal{H}}

\newcommand{\cV}{\mathcal{V}}
\newcommand{\cW}{\mathcal{W}}
\newcommand{\cI}{\mathcal{I}}

\newcommand{\cG}{\mathcal{G}}
\newcommand{\cX}{\mathcal{X}}

\newcommand{\cA}{\mathcal{A}}
\newcommand{\cY}{\mathcal{Y}}
\newcommand{\CE}{\mathcal{CE}}

\newcommand{\cF}{\mathcal{F}}

\newcommand{\cP}{\mathscr{P}}
\newcommand{\hrho}{\hat\rho}
\newcommand{\brho}{\bar\rho}
\newcommand{\bpsi}{\bar\psi}
\newcommand{\bPsi}{\bar\Psi}

\renewcommand{\tilde}{\widetilde}

\begin{document}

\title
[{Entropy gradient flows for Markov chains}]{Gradient flows of the entropy for finite Markov chains}

\author{Jan Maas}
\address{
University of Bonn\\
Institute for Applied Mathematics\\
Endenicher Allee 60\\
53115 Bonn\\
Germany}
\email{maas@iam.uni-bonn.de \hfill http://www.janmaas.org}

\thanks{The author is supported by Rubicon subsidy 680-50-0901 of the Netherlands Organisation for Scientific Research (NWO)}



 \begin{abstract}
Let $K$ be an irreducible and reversible Markov kernel on a finite set $\cX$. We construct a metric $\cW$ on the set of probability measures on $\cX$ and show that with respect to this metric, the law of the continuous time Markov chain evolves as the gradient flow of the entropy.
This result is a discrete counterpart of the Wasserstein gradient flow interpretation of the heat flow in $\R^n$ by Jordan, Kinderlehrer, and Otto (1998).
The metric $\cW$ is similar to, but different from, the $L^2$-Wasserstein metric, and is defined via a discrete variant of the Benamou-Brenier formula. 
 \end{abstract}

\maketitle
 
\section{Introduction}

Since the seminal work of Jordan, Kinderlehrer and Otto \cite{JKO98}, it is known that the heat flow on $\R^n$ is the gradient flow of the Boltzmann-Shannon entropy with respect to the $L^2$-Wasserstein metric on the space of probability measures on $\R^n$. This discovery has been the starting point for many developments in evolution equations, probability theory and geometry. We refer to the monographs \cite{AGS08,Vil03,Vil09} for an overview.
By now a similar interpretation of the heat flow has been established in a wide variety of settings, including Riemannian manifolds \cite{Erb10}, Hilbert spaces \cite{ASZ09}, Wiener spaces \cite{FSS10}, Finsler spaces \cite{OhSt09}, Alexandrov spaces \cite{GKO10} and metric measure spaces \cite{Gi10,Sav07}. 

Let $(K(x,y))_{x,y \in \cX}$ be an irreducible and reversible Markov transition kernel on a finite set $\cX$, and consider the continuous time semigroup $(H(t))_{t\geq0}$ associated with $K$. This semigroup is defined by $H(t) = e^{t(K - I)}$,
and can be interpreted as the `heat semigroup' on $\cX$ with respect to the geometry determined by the Markov kernel $K$. Therefore it seems natural to ask whether the heat flow can also be identified as the gradient flow of an entropy functional with respect to some metric on the space of probability densities on $\cX$. Unfortunately, it is easily seen that the $L^2$-Wasserstein metric over a discrete space is not appropriate for this purpose.
In fact, since the metric derivative of the heat flow in the Wasserstein metric is typically infinite in a discrete setting, 
the heat flow can not be interpreted as the gradient flow of \emph{any} functional in the $L^2$-Wasserstein metric. 
(We refer to Section \ref{sec:two-point} for a more detailed discussion.) 

 The main contribution of this paper is the construction of a metric $\cW$ on the space of probability densities on $\cX$, which allows to extend the interpretation of the heat flow as the gradient flow of the entropy to the setting of finite Markov chains.

\subsection*{Notation}  

As before, let $K : \cX \times \cX \to \R$ be a Markov kernel on a finite space $\cX$, i.e., 
\begin{align*}
 K(x,y) \geq 0 \quad \forall x,y \in \cX\;, 
 \qquad
 \sum_{y \in \cX} K(x,y) = 1 \quad\forall x\in \cX\;. 
\end{align*}
We assume that $K$ is \emph{irreducible}, which implies the existence of a unique  steady state $\pi$.
Thus $\pi$ is a probability measure on $\cX$, represented by a row vector that is invariant under right-multiplication by $K$:
\begin{align*}
  \pi(y) = \sum_{x \in \cX} \pi(x) K(x,y)\;.
\end{align*}
It follows from elementary Markov chain theory that $\pi$ is strictly positive.
We shall assume that $K$ is \emph{reversible}, i.e., $\pi(x)K(x,y) = \pi(y)K(y,x)$ for any $x, y\in \cX$. 
Consider the set
\begin{align*}
 \cP(\cX) := \Big\{ \, \rho : \cX \to \R \ | \ 
    \rho(x) \geq 0 \quad \forall x \in \cX\   \;; \ \sum_{x \in \cX} \pi(x) \rho(x)  = 1 \, \Big\}
\end{align*}
consisting of all \emph{probability densities} on $\cX$. The subset consisting of those probability densities that are strictly positive is denoted by $\cP_*(\cX)$. The \emph{relative entropy} of a probability density $\rho \in \cP(\cX)$ with respect to $\pi$ is defined by
\begin{align} \label{eq:entropy}
 \cH(\rho) = \sum_{x \in \cX} \pi(x) \rho(x) \log \rho(x)\;.
\end{align}
with the usual convention that $\rho(x) \log \rho(x) = 0$ if $\rho(x) = 0$.

\subsection*{Wasserstein-like metrics in a discrete setting}

To motivate the definition of the metric $\cW$, recall that for probability densities $\rho_0, \rho_1$ on $\R^n$, the Benamou-Brenier formula \cite{BB00} asserts that the squared Wasserstein distance $W_2$ satisfies the identity
\begin{align}\label{eq:Benamou-Brenier}
 W_2(\rho_0, \rho_1)^2
   = \inf_{\rho, \psi}
    \bigg\{ 
    \int_0^1 \int_{\R^n} |\nabla \psi_t(x)|^2 \, \rho_t(x) \dd x \dd t
    \bigg\}\;,
\end{align}
where the infimum runs over sufficiently regular curves $\rho : [0,1] \to \cP(\R^n)$ and $\psi : [0,1] \times \R^n \to \R$ satisfying the continuity equation 
\begin{equation}\begin{aligned}\label{eq:continuityEq}
 \left\{ \begin{array}{l}
 \partial_t \rho + \nabla \cdot (\rho \nabla \psi) = 0\;,\\
\rho(0) = \rho_0\;, \quad \rho(1)  = \rho_1\;.\end{array} \right.
\end{aligned}\end{equation}
Here, by a slight abuse of notation, $\cP(\R^n)$ denotes the set of probability densities on $\R^n$.
At least formally, the Benamou-Brenier formula has been interpreted by Otto \cite{O01} as a Riemannian metric on the space of probability densities on $\R^n$.

In the discrete setting, we shall define a class of pseudo-metrics $\cW$ (i.e., metrics which possibly attain the value $+\infty$) by mimicking the formulas \eqref{eq:Benamou-Brenier} and \eqref{eq:continuityEq}.

 In order to obtain a metric with the desired properties, it turns out to be necessary to define, for $\rho \in \cP(\cX)$ and $x,y \in \cX$,
\begin{align*}
 \rho(x,y) := \theta(\rho(x), \rho(y))\;,
\end{align*}
where $\theta : \R_+ \times \R_+ \to \R_+$ is a function satisfying  (A1) -- (A7) below. 
At this stage we remark that typical examples of admissible functions are the logarithmic mean $\theta(s,t) = \int_0^1 s^{1-p} t^p \dd p$, the geometric mean $\theta(s,t) = \sqrt{st}$ and, more generally, the functions $\theta(s,t) = s^\alpha t^\alpha$ for $\alpha > 0$.

Now we are ready to state the definition of $\cW$:

\begin{definition*}\label{def:new-metric}
For $\rho_0, \rho_1 \in \cP(\cX)$ we set
\begin{align*}
 \cW(\rho_0, \rho_1)^2
   := \inf_{\rho, \psi} 
   \bigg\{  \frac12   \int_0^1 
  \sum_{x,y\in \cX} (\psi_t(x) - \psi_t(y))^2
    		 K(x,y) \rho_t(x,y) \pi(x)
      \dd t 
          \bigg\}\;,
\end{align*}
where the infimum runs over all piecewise $C^1$ curves $\rho : [0,1] \to \cP(\cX)$ and all measurable functions $\psi : [0,1] \to \R^\cX$ satisfying, for a.e. $t \in [0,1]$,
\begin{align} \label{eq:cont}
 \begin{cases}
 \displaystyle\ddt \rho_t(x) 
   + \displaystyle\sum_{y \in \cX} ( \psi_t(y) - \psi_t(x) )K(x,y) \rho_t(x,y) = 0\qquad \forall x \in \cX\;, \\ 
  \rho(0) = \rho_0\;, \qquad \rho(1)  = \rho_1\;.
 \end{cases}
\end{align}
\end{definition*}

\begin{remark*}\label{rem:non-local}
Similar to the Wasserstein metric, $\cW(\rho_0, \rho_1)^2$ can be interpreted as the cost of transporting mass from its initial configuration $\rho_0$ to the final configuration $\rho_1$. However, unlike the Wasserstein metric, the cost of transporting a unit mass from $x$ to $y$ depends on the amount of mass already present at $x$ and $y$. In a continuous setting, metrics with these properties have been studied in the recent papers \cite{CLSS10,DNS09}. The essential new feature of the metric considered in this paper is the fact that the dependence is \emph{non-local}.
\end{remark*}

In order to state the first main result of the paper, we introduce some notation. Fix a probability density $\rho \in \cP(\cX)$. We shall write $x \sim_\rho y$ if $x, y \in \cX$ belong to the same connected component of the support of $\rho$. More formally, we say that $x \sim_\rho y$ if $x = y$, or if there exist $k \geq 1$ and $x_1, \ldots, x_k \in \cX$ such that
\begin{align*}
  \rho(x,x_1)  K(x,x_1), 
 \rho(x_1,x_2)K(x_1,x_2),
\;\ldots\;,
\rho(x_k,y)K(x_k,y) > 0\;.
\end{align*}

Furthermore, we set
\begin{align*}
C_\theta := \int_0^1 \frac{1}{\sqrt{\theta(1-r, 1+r)}} \dd  r \in [0, \infty] \;.
\end{align*}
It turns out that $C_\theta$ is the $\cW$-distance between a Dirac mass and the uniform density on a two-point space $\{a,b\}$ endowed with the Markov kernel defined by $K(a,b) = K(b,a) =\frac12$. Note that $C_\theta$ is finite if $\theta$ is the logarithmic or geometric mean. If $\theta(s,t) = s^\alpha t^\alpha$, then $C_\theta$ is finite for $0 < \alpha < 2$ and infinite for $\alpha \geq 2$.

For $\sigma \in \cP(\cX)$ we shall write
\begin{align*}
\cP_\sigma(\cX) := \{ \rho \in \cP(\cX) \ : \
   \cW(\rho, \sigma) < \infty \}\;.
\end{align*}

The first main result of this paper reads as follows:

\begin{theorem}\label{thm:main-metric}
The following assertions hold:
\begin{enumerate}
\item
 $\cW$ defines a pseudo-metric on $\cP(\cX)$.
\item 
\begin{itemize}
\item If $C_\theta < \infty$, then $\cW(\rho_0, \rho_1) < \infty$ for all $
\rho_0, \rho_1 \in \cP(\cX)$.
\item If $C_\theta = \infty$, the following are equivalent for $\rho_0, \rho_1 \in \cP(\cX)$: 
 \begin{enumerate}
 \item $\cW(\rho_0, \rho_1) < \infty$\;;
 \item For all $x\in \cX$ we have
\begin{align*}
 \sum_{y\sim_{\rho_0} x} \rho_0(y) \pi(y) = 
 \sum_{y\sim_{\rho_1} x} \rho_1(y) \pi(y)\;.
\end{align*}
 \end{enumerate}
\end{itemize}
\item For all $\sigma \in \cP(\cX)$, $\cW$ metrizes the topology of weak convergence on $\cP_\sigma(\cX)$.  
\item 
\begin{itemize}

\item If $C_\theta < \infty$ and $\theta$ is concave, the metric space $(\cP_*(\cX), \cW)$ is a Riemannian manifold.
\item If $C_\theta = \infty$, the metric space $(\cP_\sigma(\cX), \cW)$ is a complete Riemannian manifold for all $\sigma \in \cP(\cX)$.

\end{itemize}
\end{enumerate}
\end{theorem}

\begin{remark*}[Finiteness]\label{rem:finiteness}

Part (2) of the theorem above provides a complete characterisation of finiteness of $\cW$ for general Markov kernels, in terms of the behaviour of $\cW$ for kernels on a two-point space.
If $C_\theta = \infty$, the statement can be rephrased informally by saying that the distance $\cW(\rho_0, \rho_1)$ is finite if and only if the following conditions hold: $\rho_0$ and $\rho_1$ have equal support, and both measures assign the same mass to each connected component of their support.
In particular, it is important to note that the distance between two \emph{strictly} positive densities is finite.
\end{remark*}

\begin{remark*}[Weak convergence]\label{rem:weak}
Although (3) asserts that $\cW$ metrizes the topology of weak convergence on $\cP_\sigma(\cX)$ for every $\sigma \in \cP(\cX)$, it follows from (2) that $\cW$ does \emph{not} metrize this topology on the full space $\cP(\cX)$ if $C_\theta = \infty$. In fact, a weakly convergent sequence in $\cP_\sigma(\cX)$ converges in $\cW$-metric if and only if the weak limit belongs to $\cP_\sigma(\cX)$.
\end{remark*}

\begin{remark*}[Non-compactness]\label{rem:compact}
If $C_\theta = \infty$, we hasten to point out that the Riemannian manifold $(\cW, \cP_\sigma(\cX))$ can be a singleton. According to (2), this happens if and only if $K(x,y) \sigma(x,y) = 0$ for every $x \in \supp \sigma$ and every $y \in \cX$, which is for instance the case if $\sigma$ is the density of a Dirac measure.
If $\cP_\sigma(\cX)$ consists of more than one element, it turns out that $(\cP_\sigma(\cX), \cW)$ is non-compact. By contrast, the $L^2$-Wasserstein space over a compact metric space is compact.
\end{remark*}

\begin{remark*}[Riemannian metric]
The Riemannian metric on $(\cP_*(\cX), \cW)$ is a natural discrete analogue of the formal Riemannian metric on the Wasserstein space over $\R^n$. 
In fact, consider a smooth curve $ (\rho_t)_{t \in [0,1]}$ in $\cP_*(\cX)$ and take $t \in [0,1]$. In Section \ref{sec:metric} we shall prove that there exists a unique discrete gradient $\nabla \psi_t = (\psi_t(x) - \psi_t(y))_{x,y \in \R^n}$ such that the continuity equation \eqref{eq:cont} holds.
In view of this observation, we shall identify the tangent space at $\rho \in \cP_*(\cX)$ with the collection of discrete gradients 
\begin{align*}
 T_\rho := \{ \nabla \psi \in \R^{\cX \times \cX} \ : \ \psi \in \R^\cX \}\;.
\end{align*}
We shall regard the discrete gradient $\nabla \psi_t$ as being the tangent vector along the curve $t \mapsto \rho_t$.
The distance $\cW$ is the Riemannian distance induced by the inner product $\ip{\cdot, \cdot}_\rho$ on $T_\rho$ given by
\begin{align*}
 \ip{\nabla \phi, \nabla \psi}_\rho 
   = \frac12 \sum_{x,y\in \cX} 
 (\phi(x) - \phi(y))(\psi(x) - \psi(y)) K(x,y) \rho(x,y)\pi(x)\;.
\end{align*}
This formula is analogous to the corresponding expression in the continuous case \cite{O01}.
In Section \ref{sec:metric} we obtain a similar description of the Riemannian metric on each of the components of $\cP(\cX)$. If $\rho$ is not strictly positive, the tangent space shall be identified with the collection of discrete gradients of an appropriate subset of functions on $\cX$.

\end{remark*}

\begin{remark*}[Two-point space]\label{rem:two-point}
If $K$ is a reversible Markov kernel on a space $\cX$ consists of only two points, it is possible to obtain an explicit formula for the metric $\cW$. We refer to Section \ref{sec:two-point} for an extensive discussion. 
\end{remark*}

\begin{example*}
If $C_\theta = \infty$, it follows from Theorem \ref{thm:main-metric} that the incidence graph associated with the Markov kernel $K$ determines the topology of $(\cP(\cX),\cW)$. Let us illustrate this fact by two simple examples on a three-point space $\cX = \{ x_1, x_2, x_3 \}$. 

If $K(x_i, x_j) > 0$ for all $i \neq j$, then the space $\cP(\cX)$ consists of $7$ distinct Riemannian manifolds:
\begin{itemize}
\item one 2-dimensional manifold: $\cP_*(\cX)$;
\item three 1-dimensional manifolds: for $i = 1,2,3$,
\begin{align*}
C_i:=  \{ \rho  \in \cP(\cX) \ : \ \rho(x_j) = 0 \text{ iff } j = i \}\;.
\end{align*}
\item three singletons: for $i = 1,2,3$, 
\begin{align*}
D_i := \{ \rho  \in \cP(\cX) \ : \ \rho(x_j) = 0 \text{ iff } j \neq i \}\;.
\end{align*}
\end{itemize}

If $K(x_1, x_2), K(x_2, x_3) > 0$ and $K(x_1, x_3) = 0$, then the space $\cP(\cX)$ consists of infinitely many distinct Riemannian manifolds:
\begin{itemize}
\item one 2-dimensional manifold: $\cP_*(\cX)$;
\item two 1-dimensional manifolds: $C_1$ and $C_3$;
\item infinitely many singletons: the 
three singletons $D_i$ for $i = 1,2,3$, and the infinite collection
\begin{align*}
 \{  \{ \rho \} \ : \ \rho(x_1) > 0, \; \rho(x_3) > 0, \;\rho(x_2) = 0 \}\;.
\end{align*}
\end{itemize}
\end{example*}

\subsection*{The gradient flow of the entropy}

Since the entropy functional $\cH$ restricts to a smooth functional on the Riemannian manifold $(\cP_*(\cX), \cW)$, it makes sense to consider the associated gradient flow. Let $D_t \rho$ denote the tangent vector field along a smooth curve $\rho : (0,\infty) \to \cP_*(\cX)$ and let $\grad \phi$ denote the gradient of a smooth functional $\phi : \cP_*(\cX) \to \R$.

Consider the continuous time Markov semigroup $H(t) = e^{t(K-I)}$, $t \geq 0$, associated with $K$. It follows from the theory of Markov chains that $H(t)$ maps $\cP(\cX)$ into $\cP_*(\cX)$.
The second main result of this paper asserts that the `heat flow' determined by $H(t)$ is the gradient flow of the entropy $\cH$ with respect to $\cW$, if $\theta$ is the logarithmic mean.

\begin{theorem}[Heat flow is gradient flow of entropy]
	\label{thm:main-gradflow}
Let $\theta$ be the logarithmic mean.
For $\rho \in \cP(\cX)$ and $t \geq 0$, set $\rho_t := e^{t(K - I)} \rho$.
Then the gradient flow equation
\begin{align*}
  D_t \rho = - \grad \cH(\rho_t)
\end{align*}
holds for all $t > 0$.
\end{theorem}

\begin{remark*}
The choice of the logarithmic mean is essential in Theorem \ref{thm:main-gradflow} if one wishes to identify the heat flow as the gradient flow of the entropy associated with the function $f(\rho) = \rho \log \rho$. In section \ref{sec:gradFlow} we prove that analogous results can be proved for certain different functions $f$, if one replaces the logarithmic mean by $\theta(s,t) = \frac{s-t}{f'(s) - f'(t)}$.
The appearance of the logarithmic mean in discrete heat flow problems is not surprising. In fact, the ``Log Mean Temperature Difference'', usually called LMTD, plays an important r\^ole in the engineering literature on heat and mass transfer problems (see, e.g., \cite{McA54}), in particular in heat flow through long cylinders (see also \cite[Section 4.5]{Bha07} for a discussion). 
\end{remark*}

\begin{remark*}
For Markov chains on a two-point space $\{-1,1\}$ we shall show in Section \ref{sec:two-point} that (under mild additional assumptions) the metric $\cW$ is the \emph{unique} metric for which the gradient flow of the entropy coincides with the heat flow. We refer to Proposition \ref{prop:metric-uniqueness} below for a precise statement.
\end{remark*}

\subsection*{Ricci curvature in a discrete setting} 

A synthetic theory of Ricci curvature in metric measure spaces has been developed recently by Lott-Sturm-Villani \cite{LV09,S06}. These authors defined lower bounds on the Ricci curvature of a \emph{geodesic} metric measure space in terms of convexity properties of the entropy functional along geodesics in the $L^2$-Wasserstein metric. 
For long there has been interest to define and study a notion of Ricci curvature on discrete spaces, but unfortunately the Lott-Sturm-Villani definition cannot be applied directly. The reason is that geodesics in the $L^2$-Wasserstein space do typically not exist if the underlying metric space is discrete, even in the simplest possible example of the two-point space (see Section \ref{sec:two-point} below for more details).

The metric $\cW$ constructed in this paper does not have this defect. By a lower-semicontinuity argument it can be shown that every pair of probability densities in $\cP(\cX)$ can be joined by a constant speed geodesic.
Since $\cW$ takes over the r\^ole of the $L^2$-Wasserstein metric if $\theta$ is the logarithmic mean, the following modification of the Lott-Sturm-Villani definition of Ricci curvature seems natural:

\begin{definition}[Ricci curvature lower bound]\label{def:Ricci}
Let $K = (K(x,y))_{x, y \in \cX}$ be an irreducible and reversible Markov kernel on a finite space $\cX$. Then $K$ is said to have \emph{Ricci curvature bounded from below by $\kappa \in \R$}, if for every  $\bar\rho_0, \bar\rho_1 \in \cP(\cX)$ there exists a constant speed geodesic $(\rho_t)_{t \in [0,1]}$ in $(\cP(\cX), \cW)$ satisfying $\rho_0 = \bar\rho_0$, $\rho_1 = \bar\rho_1$, and
\begin{align*}
  \cH(\rho_t) \leq
  (1-t) \cH(\rho_0) + t \cH(\rho_1) - \frac\kappa{2} t(1-t) \cW(\rho_0, \rho_1)^2
\end{align*}
for all $t \in [0,1]$.
We set 
\begin{align*}
 \Ric(K) := \sup \{\kappa \in \R \ : \ K \text{ has Ricci curvature bounded from below by $\kappa$.} \}
\end{align*}
\end{definition}

Calculating or estimating $\Ric(K)$ in concrete situations does not appear to  be an easy task. We shall address this topic in a forthcoming publication.

Several other approaches to Ricci curvature in a discrete setting have been considered recently.

Bonciocat and Sturm \cite{BS09} adapted the definition based on displacement convexity of the entropy from \cite{LV09,S06} to the discrete setting. The non-existence of geodesics in the $L^2$-Wasserstein space is circumvented by considering approximate midpoints between measures in the $L^2$-Wasserstein metric. Using this approach it is shown that certain planar graphs have non-negative Ricci curvature.

Ollivier \cite{Oll07,Oll09} defined a notion of Ricci curvature by comparing transportation distances between small balls and their centers. This notion coincides with the usual notion of Ricci curvature lower boundedness on Riemannian manifolds and is very well adapted to study Ricci curvature on discrete spaces. In particular, it is easy to show that the Ricci curvature of the $n$-dimensional discrete hypercube is proportional to $\frac{1}{n}$. However, as has been discussed in \cite{OV10}, the relation with displacement convexity remains to be clarified.

Very recently Y. Lin and S.-T. Yau \cite{LY10} studied Ricci curvature on graphs by taking a characterisation in terms of the heat semigroup due to Bakry and Emery as a definition. With this definition it is shown that the Ricci curvature on locally finite graphs is bounded from below by $-1$.

\subsection*{Structure of the paper} 
Section \ref{sec:two-point} contains a detailed analysis of the metric $\cW$ associated with Markov kernels on a two-point space. In section \ref{sec:metric} we study the metric $\cW$ in a general setting and prove Theorem \ref{thm:main-metric}. 
In Section \ref{sec:gradFlow} we study gradient flows and present the proof of Theorem \ref{thm:main-gradflow}.

\subsubsection*{Note added}

After completion of this paper, the author has been informed about the recent preprint \cite{CHLZ11} where a related class of a metrics has been studied independently. The results obtained in both papers are largely complementary.

\subsection*{Acknowledgement}
{\small
The author is grateful to Matthias Erbar, Nicola Gigli, Nicolas Juillet, Giuseppe Savar\'e, and Karl-Theodor Sturm for stimulating discussions on this paper and related topics.} 

\section{Analysis on the two-point space} 
\label{sec:two-point}

In this section we shall carry out a detailed analysis of the metric $\cW$ in the simplest case of interest, where the underlying space is a two-point space, say $\cX = \cQ^1 =  \{ a, b\}$. 
The reason for discussing the two-point space separately is twofold. Firstly,  it is possible to perform explicit calculations, which lead to simple proofs and more precise results than in the general case. Secondly, some of the results obtained in this section shall be used in Section \ref{sec:metric}, where results for more general Markov chains are obtained by comparison arguments involving Markov chains on a two-point space. 

\subsection*{Markov chains on the two-point space}
Consider a Markov kernel $K$ with transition probabilities 
\begin{align} \label{eq:K}
 K(a, b) = p\;, \qquad
 K(b, a) = q\;
\end{align}
for some $p, q \in (0,1]$. Then the associated continuous time semigroup  $H(t) = e^{t(K-I)}$ is given by
\begin{align*}
 H(t) 
   = \frac{1}{p+q} \bigg( 
    \left[\begin{array}{cc}q & p \\q & p\end{array}\right]
    +  e^{-(p+q)t} 
    \left[\begin{array}{cc}p & -p \\-q & q\end{array}\right] 
    \bigg)\;.
\end{align*}
and the stationary distribution $\pi$ satisfies 
\begin{align*}
  \pi(a) =  \frac{q}{p+q}\;,
  \qquad 
   \pi(b) =  \frac{p}{p+q}\;.
\end{align*}
Since $K(a,b) \pi(a) = K(b,a) \pi(b)$, we observe that $K$ is reversible.
Every probability measure on $ \cQ^1$ is of the form $\frac12 ((1- \beta) \delta_{a} + (1+\beta)\delta_{b})$ for some $\beta \in [-1,1]$. The corresponding density $\rho^\beta$ with respect to $\pi$ is then given by
\begin{align*}
  \rho^\beta(a) :=  \frac{p + q}{q} \frac{1- \beta}{2}  \;,\qquad
   \rho^\beta(b) :=  \frac{p + q}{p} \frac{1+ \beta}{2} \;.
\end{align*}
It follows that $H(t) \rho^\beta = \rho^{\beta_t}$ where
\begin{align} \label{eq:heat-flow}
\beta_t :=
 \frac{p-q}{p+q}\big(1 - e^{-(p+q)t}\big) + \beta e^{-(p+q)t}\;,
\end{align}
thus $\beta$ solves the differential equation
\begin{align} \label{eq:beta-ode}
  \dot\beta_t =  p(1- \beta_t) -  q(1 + \beta_t)\;.
\end{align}

\begin{remark}[Limitations of the $L^2$-Wasserstein distance]\label{rem:Wasserstein}
Before introducing a new class of (pseudo)-metrics on $\cP(\cQ^1)$, we shall argue why the $L^2$-Wasserstein metric $W_2$ is not appropriate for the purposes of this paper.
First we shall show that -- as we already mentioned in the introduction -- the metric derivative of the heat flow is infinite with respect to the $L^2$-Wasserstein metric. To see this, take $\beta \in [-1,1] \setminus \{ \frac{p-q}{p+q}\}$, and let $u(t) := H(t) \rho^\beta = \rho^{\beta_t}$ be the heat flow starting at $\rho^\beta$. Since $W_2(\rho^\alpha, \rho^\beta) = \sqrt{2|\beta - \alpha|}$ for $\alpha, \beta \in [-1,1]$, we have
\begin{align*}
  |\dot u|(t) & := \limsup_{s \to t} \frac{ W_2( u(t), u(s))}{|t - s|}
          =  \sqrt{2} \limsup_{s\to t} 
          	\frac{\sqrt{|\beta_t - \beta_s|}}{|t - s|} 
	\\&	  =  \sqrt{2 \Big| \beta - \frac{p-q}{p+q}\Big| } \limsup_{s\to t} 
		  	\frac{\sqrt{|
					e^{-(p+q)t} - e^{-(p+q)s} | }}{|t - s|} 
          = + \infty\;.
\end{align*}
In particular, the heat flow is \emph{not} a curve of maximal slope (see, e.g., \cite{AGS08} for this concept of gradient flow) for any functional on $\cP(\cQ^1)$.

Furthermore, the Lott-Sturm-Villani definition of Ricci curvature \cite{LV09,S06} cannot be applied in the discrete setting, since $W_2$-geodesics between distinct elements of $\cP(\cQ_1)$ do not exist. To see this, let $\{ \rho^{\beta(t)}\}_{0 \leq t \leq 1}$ be a constant speed geodesic in $\cP(\cQ^1)$. For $s,t\in [0,1]$ we then have
\begin{align*}
 \sqrt{2|\beta(t) - \beta(s)|} 
 &= W_2(\rho^{\beta(t)}, \rho^{\beta(s)}) 
 \\&= |t-s| W_2(\rho^{\beta(0)}, \rho^{\beta(1)}) 
 = |t-s| \sqrt{2|\beta(0) - \beta(1)|}\;,
\end{align*}
which implies that $t \mapsto \beta(t)$ is $2$-H\"older, hence constant on $[0,1]$. It thus follows that all constant speed $W_2$-geodesics are constant.
\end{remark}

\subsection*{A new metric}
Given a fixed Markov chain $K$ on $\{a,b\}$ we shall define a (pseudo-)metric $\cW$ on $\cP(\{a,b\})$ that depends on the choice of a function $\theta : \R_+ \times \R_+ \to \R_+$. 
The following assumptions will be in force throughout this section:

\begin{assumption} \label{ass:theta}
The function $\theta : [0,\infty) \times [0,\infty) \to [0,\infty)$ has the following properties:
\begin{itemize}
\item[(A1)] $\theta$ is continuous on $[0,\infty) \times [0,\infty)$;
\item[(A2)] $\theta$ is continuously differentiable on $(0,\infty) \times (0,\infty)$;
\item[(A3)] $\theta(s,t) = \theta(t,s)$ for $s, t \geq 0$;
\item[(A4)] $\theta(s,t) > 0$ for $s,t > 0$.
\end{itemize}
\end{assumption}

The most interesting choice for the purposes of this paper is the case where $\theta$ is the logarithmic mean defined by $\theta(s, t) := \int_0^1 s^{1-p} t^p \dd p$.

To simplify notation we define, for $\beta \in [-1,1]$,
\begin{align*}
 \hrho(\beta) = \theta(\rho^\beta(a), \rho^\beta(b))\;.
\end{align*}
On the two-point space the variational definition of $\cW$ given in the introduction can be simplified as follows:

\begin{lemma}\label{lem:var-form}
For $\alpha, \beta \in [-1,1]$ we have
\begin{align} \label{eq:var-2}
 \cW(\rho^\alpha, \rho^\beta)^2
   = \inf_{\beta} \bigg\{ \frac{p+q}{4pq}  
      \int_0^1  \frac{\dot\beta_t^2}{\hrho(\beta_t)} \one_{\{\hrho(\beta_t) > 0\}} \dd t \bigg\}\;,
\end{align}
where the infimum runs over all piecewise $C^1$-functions $\beta: [0,1] \to [-1,1]$.
\end{lemma}

\begin{proof}
Substituting $\chi(t) = \psi_t(b) - \psi_t(a)$ in the definition of $\cW$, one obtains
\begin{align*}
 \cW(\rho^\alpha, \rho^\beta)^2
   = \inf_{\beta,\chi} \bigg\{ \frac{pq}{p+q} \int_0^1 \hrho(\beta_t) \chi_t^2  \dd t \bigg\}\;,
\end{align*}
where the infimum runs over all piecewise $C^1$-functions $\beta: [0,1] \to [-1,1]$ and all measurable functions $\chi : [0,1] \to \R$ satisfying $\beta_0 = \alpha$, $\beta_1 = \beta$ and
\begin{align*}
 \dot\beta_t = \frac{2pq}{p+q} \hrho(\beta_t) \chi_t\;.
\end{align*}
The result follows by inserting the latter constraint in the expression for $ \cW(\rho^\alpha, \rho^\beta)$. 
\end{proof}

Lemma \ref{lem:var-form} provides a representation of $\cW(\rho^\alpha, \rho^\beta)$ in terms of a one-dimensional variational problem. Note that some care needs to be taken when solving this problem, since for some choices of $\theta$ (including the logarithmic mean) the denominator in \eqref{eq:var-2} tends to $0$ as $\beta_t$ tends to $\pm 1$.
The following result provides an explicit formula for $\cW$:

\begin{theorem}\label{thm:two-point}
For $-1 \leq \alpha \leq \beta \leq 1$ we have
\begin{align*}
 \cW(\rho^\alpha, \rho^\beta) 
   =  \frac12 \sqrt{ \frac1p + \frac1q } 
  \int_\alpha^\beta \frac{1}{\sqrt{\hrho(r)}}
  \dd  r \in [0, \infty]\;.
\end{align*}
\end{theorem}

\begin{proof}
Suppose first that $\alpha$ and $\beta$ belong to $(-1,1)$. (If $\hrho$ is bounded away from $0$, this distinction is not necessary.)
It is easily checked that the infimum in \eqref{eq:var-2} may be restricted to monotone functions $\beta$. 
Since $g : r \mapsto \frac{1}{\hrho(r)}$ is bounded on compact intervals in $(-1,1)$, \eqref{eq:var-2} reduces to an elementary one-dimensional variational problem, which admits a minimizer, say $\xi$, that solves the Euler-Lagrange equation
\begin{align*}
  2 \ddot\xi_t g(\xi_t) + \dot \xi_t ^2 g'(\xi_t) = 0\;.
\end{align*}
This equation implies that $t \mapsto \dot \xi_t \sqrt{g(\xi_t)}$ is constant, say equal to $C$. Since $\alpha \leq \beta$, it follows that $C > 0$. We infer that
\begin{align*}
 \cW(\rho^\alpha, \rho^\beta)^2
   = \frac{p+q}{4pq} \int_0^1 \frac{\dot\xi_t^2}{\hrho(\xi_t)}  \dd t
   =  \frac{p+q}{4pq}C^2\;.
\end{align*}
Moreover, $\xi$ is monotone, hence invertible. It follows from the inverse function theorem that its inverse $\gamma : [\alpha, \beta] \to [0,1]$ satisfies $\gamma'(r) = C^{-1} \sqrt{g(r)}$. We thus obtain
\begin{align*}
 1 = \gamma(\beta) - \gamma(\alpha)
   = \int_\alpha^\beta \gamma'(r) \dd r
   = C^{-1} \int_\alpha^\beta \sqrt{g(r)} \dd r\;,
\end{align*}
hence 
\begin{align*}
 \cW(\rho^\alpha, \rho^\beta)
   = \frac{C}{2} \sqrt{\frac1p + \frac1q}
   = \frac{1}{2} \sqrt{\frac1p + \frac1q}
   		 \int_\alpha^\beta \sqrt{g(r)} \dd r\;,
\end{align*}
which implies the desired identity.

The general case $-1 \leq \alpha \leq \beta \leq 1$ follows from a straightforward continuity argument.
\end{proof}

For $\beta \in [-1,1]$ it will be useful to define
\begin{align} \label{eq:def-phi}
 \phi(\beta) := \frac12 \sqrt{ \frac1p + \frac1q } 
  \int_0^\beta \frac{1}{\sqrt{\hrho(r)}}
  \dd  r \in [-\infty, \infty] \;,
\end{align}
so that Theorem \ref{thm:two-point} implies that 
\begin{align*}
 \cW(\rho^\alpha, \rho^\beta) = |\phi(\alpha) - \phi(\beta)|
\end{align*}
for $\alpha, \beta \in [-1,1]$.
It follows from the assumption on $\theta$ that $\phi$ is real-valued, continuous and strictly increasing on $(-1,1)$. Moreover, $\phi(\pm 1) = \lim_{\beta \to \pm 1} \phi(\beta)$ is possibly $\pm \infty$, depending on the behaviour of $\theta$ near $0$.

In order to avoid having to distinguish between several cases in the results below, we set
\begin{align*}
(-1,1)_* = \{ \beta \in [-1,1] : |\phi(\beta)|< \infty \}\;, 
\qquad I = \{ \phi(\beta) : \beta \in (-1,1)_* \}\;,
\end{align*}
and
\begin{align*}
 \cP_{\one}(\cQ^1) := \{ \rho^\beta \in \cP(\cQ^1) : \beta \in (-1,1)_* \}\;.
\end{align*}

It follows from the remarks above that $(-1,1)\subseteq (-1,1)_* \subseteq [-1,1]$ and that $I$ is a (possibly infinite) closed interval in $\R$.
The following result, which summarises this discussion, is now obvious:

\begin{proposition}\label{prop:metric-2pt}
The function $\cW$ defines a pseudo-metric on $\cP(\cQ^1)$ that restricts   to a metric on $\cP_{\one}(\cQ^1)$. 
The mapping
\begin{align*}
  J : \rho^\beta \mapsto \phi(\beta)
\end{align*}
defines an isometry from $(\cP_{\one}(\cQ^1), \cW)$ onto $I$ endowed with the euclidean metric. In particular,  $(\cP_{\one}(\cQ^1), \cW)$ is complete.
\end{proposition}

The most interesting case for the purposes of this paper is the following:

\begin{example}[Logarithmic mean]\label{ex:log}
If $\theta$ is the logarithmic mean, i.e., $\theta(s,t) = \int_0^1 s^{1-r} t^r  \dd r$, then $\hrho(-1) = \hrho(1) = 0$ and for $\beta \in (-1,1)$ we have
\begin{align*}
 \hrho(\beta)
  =\frac{p+q}{2pq}  
   \frac{q(1+ \beta) - p(1-\beta)}{\log q(1+\beta) - \log p(1-\beta)} \;.
\end{align*} 
In this case we have $(-1,1)_* = [-1,1]$ and $I = [\phi(-1), \phi(1)]$ is a compact interval.
Furthermore, for $-1 \leq\alpha \leq \beta \leq 1$,
\begin{align*}
 \cW(\rho^\alpha, \rho^\beta) 
  = \frac1{\sqrt2} \int_\alpha^\beta 
     \sqrt{\frac{\log q(1+r) - \log p(1-r)}{q(1+ r) - p(1-r)}} \dd r\;.
\end{align*}
If moreover $p = q$, we have
\begin{align*}
 \hrho(\beta) = 
   \frac{\beta}{\arctanh \beta}\;.
\end{align*}
and 
\begin{align*}
 \cW(\rho^\alpha, \rho^\beta) 
  = \frac1{\sqrt{2p}}\int_\alpha^\beta 
     \sqrt{\frac{\arctanh r}{r}} \dd r\;.
\end{align*}

\end{example}

Recall that a \emph{constant speed geodesic} in a metric space $(M,d)$ is a curve $u : [0,1] \to M$ satisfying 
\begin{align*}
d( u(s), u(t) )  = |t-s| d(u(0), u(1))
\end{align*}
for all $s, t \in [0,1]$.

The next result gives a characterisation of $\cW$-geodesics in $\cP_{\one}(\cQ^1)$.

\begin{proposition}[Characterisation of geodesics]\label{prop:geodesic}
Let $\rho, \sigma \in \cP_{\one}(\cQ^1)$. There exists a unique constant speed geodesic $\{ \rho^{\gamma(t)} \}_{0 \leq t \leq 1}$ in $\cP_{\one}(\cQ^1)$ with $\rho^{\gamma(0)} = \rho$ and $\rho^{\gamma(1)} = \sigma$. Moreover, the function $\gamma$ belongs to $C^1([0,1];\R)$ and satisfies the differential equation
\begin{align}\label{eq:geod}
 \gamma'(t) = 2 w        \sqrt{
 \frac{pq}{p+q}
\hrho(\gamma(t))          }
\end{align}
for $t \in [0,1]$, where $w := \sgn(\beta - \alpha) \cW(\rho^\alpha, \rho^\beta)$.
\end{proposition}

\begin{proof} 

Since the mapping $J$ is an isometry from $\cP_{\one}(\cQ^1)$ onto $I$, existence and uniqueness of geodesics follow directly from the corresponding facts in $I$.

Take now $\alpha, \beta\in (-1,1)_*$ and let $\gamma \in C^1([0,1];\R)$ be the solution to \eqref{eq:geod} with initial condition $\gamma(0) = \alpha$. For $0 \leq s < t \leq 1$ we then obtain by \eqref{eq:def-phi},
\begin{align*}
 \phi(\gamma(t)) - \phi(\gamma(s))
   = \int_s^t \phi'(\gamma(r)) \gamma'(r) \dd r
   = w(t-s)\;,
\end{align*} 
which implies that $\cW(\rho^{\gamma(t)},\rho^{\gamma(s)}) = |w|(t-s)$ and $\gamma(1) = \beta$,
hence $t \mapsto \rho^{\gamma(t)}$ is a constant speed geodesic between $\rho^\alpha$ and $\rho^\beta$. 
\end{proof}

\subsection*{Gradient flows}

In order to identify the heat flow as a gradient flow in $\cP(\cQ^1)$, we make the following assumption:
\begin{assumption}\label{ass:more}
In addition to (A1) - (A4) we assume that there exists a function $f\in C([0,\infty);\R) \cap C^2((0,\infty);\R)$ satisfying $f''(t) > 0$ for $t >0 $, and 
\begin{align}\label{eq:theta-f}
  \theta(s,t) = \frac{s-t}{f'(s) - f'(t)}, 
\end{align}
for all $s , t  >0$ with $s \neq t$.
\end{assumption}

\begin{example}\label{ex:f}
Note that this assumption is satisfied in Example \ref{ex:log} with $f(t) = t \log t$\;.
\end{example}

Consider the functional $\cF : \cP(\cQ^1) \to \R$ defined by
\begin{align*}
 \F(\rho) := \sum_{x\in \cQ^1}  f(\rho(x))\pi(x)
\end{align*}
where $f : \R_+ \to \R$ has been defined above.
It thus follows that
\begin{align} \label{eq:F-formula}
 \F(\rho^\beta) :=  \frac{q}{p+q}
    f( \rho^\beta(a) ) 
    + \frac{p}{p+q}
    f (\rho^\beta(b) ) \;.
\end{align}

Proposition \ref{prop:metric-2pt} implies that $(\cP_{\one}(\cQ^1), \cW)$ is a complete $1$-di\-men\-sional Riemannian manifold, which has a boundary if and only if $(-1,1)$ is a proper subset of $(-1,1)_*$. In particular, it makes sense to study gradient flows in $(\cP_{\one}(\cQ^1), \cW)$.

\begin{proposition}[Heat flow is the gradient flow of the entropy]\label{prop:metric}
For $\beta \in [-1,1]$ let $u : t \mapsto \rho^{\beta_t} =  H(t)\rho^\beta $ be the heat flow trajectory starting from $\rho^\beta$. Then $u$ is a gradient flow trajectory of the functional $\cF$ in the Riemannian manifold $(\cP_{\one}(\cQ^1), \cW)$.
\end{proposition}

\begin{proof}
Recall that the function $J : \rho^\beta \mapsto \phi(\beta)$
maps $\cP_{\one}(\cQ^1)$ isometrically onto a closed interval $I \subseteq \R$. 
Therefore it suffices to show that the gradient flow equation
\begin{align} \label{eq:grad-flow-1d}
 \ddt \phi(\beta_t) = - \tilde \cF'(\phi(\beta_t))
\end{align}
holds for $t > 0$, where $\tilde\cF := \cF \circ J^{-1}$.

To prove this, we set 
\begin{align*}
c_{pq} := \frac12\sqrt{ \frac1p + \frac1q}\;,\qquad
  \ell(\beta) := \rho^\beta(a) \;, \qquad 
  r(\beta)  :=\rho^\beta(b) \;,
\end{align*} 
for brevity. Using \eqref{eq:def-phi} and \eqref{eq:theta-f} we obtain
\begin{align} \label{eq:phi-prime}
\phi'(\beta) 
  = \frac{c_{pq}}{\sqrt{\hrho(\beta)}} 
  = c_{pq} \sqrt{\frac{f'(r(\beta)) - f'(\ell(\beta))}
  			{r(\beta)- \ell(\beta)}}\;.
\end{align}
Since 
\begin{align*}
 \tilde \cF(\phi(\beta)) 
 & = \tilde \cF(J(\rho^\beta)) 
  = \cF(\rho^\beta) =  \frac{q}{p+q}
    f( \ell(\beta)) 
    + \frac{p}{p+q}
    f( r(\beta))\;,
\end{align*}
it follows that $\tilde \cF$ is continuously differentiable on $I$ and 
\begin{align*}
  \tilde \cF'(\phi(\beta)) 
  & = \frac{f'(r(\beta)) - f'(\ell(\beta))}{2\phi'(\beta)}
\\&   = 
    \frac{1}{2c_{pq}}\sqrt{\big(  r(\beta)-  \ell(\beta) \big) 
    		\big(f'( r(\beta)) - f'( \ell(\beta)) \big)}\;.
\end{align*}
On the other hand, \eqref{eq:beta-ode} and \eqref{eq:phi-prime} imply that
\begin{align*}
\ddt \phi(\beta_t) 
 & = \big( p(1- \beta_t) -  q(1 + \beta_t) \big) \phi'(\beta_t)
 \\& = - \frac1{2c_{pq}^2} (r(\beta_t) - \ell(\beta_t)) \phi'(\beta_t)
 \\& = -\frac{1}{2c_{pq}}\sqrt{\big(  r(\beta_t)-  \ell(\beta_t) \big) 
    		\big(f'( r(\beta_t)) - f'( \ell(\beta_t)) \big)}\;.
\end{align*}
Combining the latter two identities we obtain \eqref{eq:grad-flow-1d}, which completes the proof.
\end{proof}

In order to investigate the convexity of $\cF$ along $\cW$-geodesics, we consider the function $K : (-1,1) \to \R$
defined by
\begin{align*}
 K(\beta) & := 
   \frac{p+q}2 + \frac12\hrho(\beta)
    \big( q f''( \rho^\beta(b)) + p f''( \rho^\beta(a)) \big)
\end{align*}
and 
\begin{align}\label{eq:kappa}
  \kappa := \inf\{K(\beta) : \beta \in (-1,1)\}\;.
\end{align}
Since $f'' > 0$, it follows that $\kappa \geq \frac{p+q}{2}$. 

\begin{remark}\label{rem:cacl}
If $f(\rho) = \rho \log \rho$, straightforward calculus shows that
\begin{align*}
  K(\beta) =
    \frac{p+q}{2} + \frac{1}{1 - \beta^2}
      \frac{ q(1+\beta) - p(1-\beta) }
    					{\log q(1+\beta) - \log p(1-\beta)}  
\end{align*}
If moreover $p =q$, one has
\begin{align*}
 K(\beta) =  
      p\bigg(1 + \frac{1}{1-\beta^2}
      	 \frac{\beta}{\arctanh \beta}\bigg)
	 \qquad \text{and} \qquad \kappa = 2p \;. 
\end{align*}
\end{remark}

It turns out that $\kappa$ determines the convexity of the functional $\cF$:

\begin{proposition}[Convexity of $\F$ along $\cW$-geodesics]\label{prop:convexity}
Let $\kappa$ be defined by \eqref{eq:kappa}.
The functional $\cF$ is $\kappa$-convex along geodesics. More explicitly, 
let $\bar\rho_0, \bar\rho_1 \in \cP_{\one}(\cQ^1)$ and let $\{\rho_t\}_{0 \leq t \leq 1}$ be the unique constant speed geodesic satisfying $\rho_0 = \bar\rho_0$ and  $\rho_1 = \bar\rho_1$. 
Then the inequality 
\begin{align*}
 \cF(\rho_t) 
   \leq (1 - t) \cF(\rho_0)  +   t \cF(\rho_1)
       - \frac{\kappa}{2} t (1-t) \cW^2(\rho_0, \rho_1)
\end{align*}
holds for all $t \in [0,1]$.
\end{proposition}

\begin{proof}
Let $\alpha, \beta \in (-1,1)_*$ be such that $\bar\rho_0 = \rho^\alpha$ and $\bar\rho_1 = \rho^\beta$ and set $w := \cW(\rho^\alpha, \rho^\beta)$. Without loss of generality we assume that $\alpha \leq \beta$. Proposition \ref{prop:geodesic} implies that $\rho_t = \rho^{\gamma(t)}$, where $\gamma$ satisfies \eqref{eq:geod}.

Set $\zeta(t) := \cF(\rho_t)$. It suffices to show that $\zeta''(t) \geq w^2\kappa $ for $t \in [0,1]$.
By \eqref{eq:F-formula} we have
\begin{align*}
 \zeta'(t) =  \frac12 \gamma'(t) 
 			\big(  f'( \rho^{\gamma(t)}(b)) - f'( \rho^{\gamma(t)}(a)) \big) \;,
\end{align*}
and therefore \eqref{eq:geod} implies that
\begin{align*}
 \zeta'(t) =   w \sqrt{\frac{pq}{ p+q }}  
 \sqrt{\big(  \rho^{\gamma(t)}(b) -  \rho^{\gamma(t)}(a)\big) 
  \big(  f'( \rho^{\gamma(t)}(b)) - f'( \rho^{\gamma(t)}(a)) \big)}\;.
\end{align*}
Differentiating this identity and using \eqref{eq:geod} once more, we obtain
\begin{align*}
 \zeta''(t) 
 =  w^2 K(\gamma(t))
 \geq  w^2 \kappa\;, 
\end{align*}
which completes the proof.
\end{proof}

The question arises whether the metric $\cW$ constructed above is the \emph{unique} geodesic metric on $\cP(\cQ^1)$ for which the heat flow is the gradient flow of the entropy. The answer is affirmative, provided that one requires that the left part $\{ \rho^\beta : \beta < \bar\beta\}$ and the right part $\{ \rho^\beta : \beta > \bar\beta\}$ of $\cP_{\one}(\cQ^1)$ are patched together in a `reasonable' way. Here $\bar\beta := \frac{p-q}{p+q}$, so that $\rho^{\bar\beta}$ corresponds to equilibrium. Such a condition is necessary, since the heat flow starting at $\rho^\beta$ with $\beta > \bar\beta$ does not `see' the measures $\rho^\alpha$ with $\alpha < \bar\beta$, and vice versa.

A precise uniqueness statement is given below. Since we shall not use this result elsewhere in the paper, we postpone its technical proof to Appendix \ref{sec:proof}, where the notions of $2$-absolute continuity and $\text{EVI}_0(\cF)$ are defined as well.

\begin{proposition}[Uniqueness of the metric]\label{prop:metric-uniqueness}
Let $\cM$ be a geodesic metric on $\cP_{\one}(\cQ^1)$ with the following properties:
\begin{enumerate}
\item For $\beta \in (-1,1)_*$, the heat flow $t \mapsto \rho^{\beta_t}$ given by \eqref{eq:heat-flow}, is a 2-absolutely continuous curve satisfying $\emph{EVI}_0(\cF)$.
\item \label{item:glueing} For $\alpha, \beta \in (-1,1)_*$ with $\alpha \leq \bar\beta \leq \beta$, we have 
\begin{align*}
\cM(\rho^{\alpha}, \rho^{\beta}) = \cM(\rho^{\alpha}, \rho^{\bar\beta}) + \cM(\rho^{\bar\beta}, \rho^\beta)\;.
\end{align*}
\end{enumerate}
Then $\cM = \cW$.
\end{proposition}

Note that (1) and (2) of Proposition \ref{prop:metric-uniqueness}  are satisfied if $\cM = \cW$. Indeed, since $\cF$ is convex by Proposition \ref{prop:convexity}, (1) follows from \cite[Proposition 23.1]{Vil09}. Furthermore (2) follows from the explicit expression of $\cW$ obtained in Theorem \ref{thm:two-point}.

\section{A Wasserstein-like metric for Markov chains}\label{sec:metric}

In this section we consider a Markov kernel $K = (K(x,y))_{x,y \in \cX}$ on a finite state space $\cX$. We assume that $K$ is irreducible, and denote its unique steady state by $\pi$. For all $x \in \cX$ we then have $\pi(x) > 0$. We also assume that $K$ is reversible, or equivalently, that the \emph{detailed balance equations}
\begin{align} \label{eq:detailed-balance}
 K(x,y) \pi(x) = K(y,x) \pi(y)
\end{align}
hold for all $x,y \in \cX$.

\subsection*{Definition of the (pseudo-)metric}

We start with the definition of a class of Wasserstein-like pseudo-metrics on $\cP(\cX)$.
As in Section \ref{sec:two-point}, the metric depends on the choice of a function $\theta : \R_+ \times \R_+ \to \R_+$, which we fix from now on. 
To simplify notation, we set
\begin{align*}
  \rho(x, y) := \theta(\rho(x), \rho(y))
\end{align*}
for $\rho \in \cP(\cX)$ and $x , y \in \cX$.

\begin{assumption}\label{ass:theta-2}
Throughout this section we shall assume that $\theta$ satisfies Assumption \ref{ass:theta}. In addition we impose the following assumptions:
\begin{itemize}
\item[(A5)] (Zero at the boundary): $\theta(0,t) = 0$ for all $t \geq 0$.
\item[(A6)] (Monotonicity): $\theta(r, t) \leq \theta(s,t)$ for all $0 \leq r \leq s$ and $t \geq 0$. 
\item[(A7)] (Doubling property): for any $T > 0$ there exists a constant $C_d > 0$ such that 
\begin{align*}
\theta(2 s, 2 t) \leq 2C_d \theta(s,t)
\end{align*}
whenever $0 \leq s,t \leq T$.
\end{itemize}
\end{assumption}

\begin{remark}\label{rem:not used}
Actually, the additional assumptions (A5) -- (A7) shall not be used until Theorem \ref{thm:finiteness}.
\end{remark}

At some places, in particular in Lemmas \ref{lem:comparison-one} and \ref{lem:comparison} below, it is possible to obtain sharper results by imposing one or both of the following assumptions as well. Note that (A7$'$) implies (A7).  

\begin{itemize}
\item[(A7$'$)] \emph{(Positive homogeneity): $\theta(\lambda s, \lambda t) = \lambda \theta(s,t)$ for $\lambda > 0$ and $s,t \geq 0$.}
\item[(A8)] \emph{(Concavity): 
the function $\theta : \R_+ \times \R_+ \to \R_+$  is concave.}
\end{itemize}
Observe that (A7$'$) and (A8) hold if $\theta$ is the logarithmic mean.

\begin{definition}[of the pseudo-metric $\cW$]\label{def:metric}
For $\bar\rho_0, \bar\rho_1 \in \cP(\cX)$ we define 
\begin{align*}
 \cW(\bar\rho_0, \bar\rho_1)^2 
  := \inf \bigg\{ &  \frac12   \int_0^1 
  \sum_{x,y\in \cX} (\psi_t(x) - \psi_t(y))^2
    		 K(x,y) \rho_t(x,y) \pi(x)
      \dd t 
         \ : \\ & \ {(\rho, \psi) \in \CE_1(\bar\rho_0, \bar\rho_1)}  \bigg\}\;,
\end{align*}
where, for $T > 0$, $\CE_T(\rho_0,\rho_1)$ denotes the collection of pairs $(\rho, \psi)$ satisfying the following conditions:
 \begin{align} \label{eq:conditions} 
 \left\{ \begin{array}{ll}
{(i)} & \rho : [0,T] \to \R^\cX  \text{ is piecewise $C^1$}\;;\\ 
{(ii)} &  \rho_0 = \bar\rho_0\;, \qquad \rho_1 = \bar\rho_1\;; \\
{(iii)} &  \rho_t \in \cP(\cX) \text{ for all $t \in [0,T]$}\;;\\
{(iv)} & \psi  : [0,T] \to \R^\cX \text{ is measurable}\;;\\
{(v)} &  \text{For all $x \in \cX$ and a.e. $t \in (0,T)$ we have}\\
  &\displaystyle{\dot \rho_t(x) 
   + \sum_{y \in \cX} ( \psi_t(y) - \psi_t(x) )K(x,y) \rho_t(x,y)
     				 = 0\;}.\
\end{array} \right.
\end{align}
\end{definition}
The latter equation may be thought of as a `continuity equation'.
For simplicity we shall often write 
\begin{align*}
\CE(\rho_0, \rho_1) := \CE_1(\rho_0, \rho_1)\;.
\end{align*} 

\begin{remark}[Matrix reformulation]\label{rem:translate}
It will be very useful to reformulate Definition \ref{def:metric} in terms of matrices.
For $\rho \in \cP(\cX)$ consider the matrices $A(\rho)$ and $B(\rho)$ in $\R^{\cX \times \cX}$ defined by
\begin{align*}
 A_{x,y}(\rho) :=  \left\{ \begin{array}{ll}
\sum_{z \neq x} K(x,z)\rho(x,z)\pi(x) \;,
 & \text{$x = y$}\;,\\
-  K(x,y)\rho(x,y) \pi(x)\;,
 & \text{$x \neq y$}\;,\end{array} \right.
\end{align*}
and
\begin{align*}
 B_{x,y}(\rho) :=  \left\{ \begin{array}{ll}
\sum_{z \neq x} K(x,z)\rho(x,z)\;,
 & \text{$x = y$}\;,\\
- K(x,y)\rho(x,y)\;,
 & \text{$x \neq y$}\;.\end{array} \right.
\end{align*}
Definition \ref{def:metric} can then be rewritten as
\begin{align}\label{eq:metric-matrix}
 \cW(\bar\rho_0, \bar\rho_1)^2 
  = \inf \bigg\{    \int_0^1 [ A(\rho_t) \psi_t, \psi_t]\dd t 
         \ : \ {(\rho, \psi) \in \CE(\bar\rho_0, \bar\rho_1)}  \bigg\}\;,
\end{align}
and the `continuity equation' in \eqref{eq:conditions} reads as
\begin{align} \label{eq:cont-equ}
 \dot\rho_t = B( \rho_t )  \psi_t\;.
\end{align}
Here and in the sequel we use square brackets $[\cdot,\cdot]$ to denote the standard inner product in $\R^\cX$.
It follows from the detailed balance equations \eqref{eq:detailed-balance} that $A(\rho)$ is symmetric, but $B(\rho)$ is not necessarily symmetric. 
Since $\sum_{y \neq x} |A_{x,y}(\rho)| = A_{x,x}(\rho) \geq 0$ for all $x \in \cX$, the matrix $A(\rho)$ is diagonally dominant, which implies that
\begin{align} \label{eq:A-nonneg}
  [A(\rho)\psi, \psi] \geq 0
\end{align} 
for all $\psi \in \R^\cX$.
Note that 
\begin{align*}
 A(\rho) = \Pi B(\rho)\;,
\end{align*}
where the diagonal matrix $\Pi \in \R^{\cX \times \cX}$ is defined by
\begin{align*}
 \Pi := \diag( \pi(x))_{x\in \cX}
\end{align*}
\end{remark}

\subsection*{Geometric interpretation}

Before continuing we present another, more geometric reformulation of Definition \ref{def:metric} which makes the connection to the Benamou-Brenier formula \ref{eq:Benamou-Brenier} (even) more apparent. We introduce some notation that will be used throughout the remainder of the paper.

For $\psi \in \R^\cX$ we consider the discrete gradient $\nabla \psi \in \R^{\cX \times \cX}$ defined by
\begin{align*}
 \nabla \psi(x,y) :=  \psi(x) - \psi(y)\;,
\end{align*}
and for $\Psi \in \R^{\cX \times \cX}$ we consider the divergence $\nabla \cdot \Psi \in \R^\cX$ defined by
\begin{align*}
( \nabla \cdot \Psi )(x) 
  := \frac12 \sum_{y \in \cX} K(x,y) (\Psi(y,x) - \Psi(x,y) ) \in \R\;.
\end{align*}
It is easily checked that the ``integration by parts formula'' holds:
\begin{align*}
 \ip{\nabla \psi, \Psi}_\pi = -\ip{\psi,\nabla\cdot \Psi}_\pi\;,
\end{align*}
where, for $\phi,\psi \in \R^\cX$ and $\Phi,\Psi \in \R^{\cX \times \cX}$,
\begin{align*}
\ip{\phi, \psi}_\pi &= \sum_{x \in \cX} \phi(x) \psi(x) \pi(x)\;, \\
\ip{\Phi, \Psi}_\pi &= \frac12 \sum_{x,y \in \cX} 
         \Phi(x,y) \Psi(x,y) K(x,y) \pi(x)\;.
\end{align*}
Furthermore, for $\rho \in \cP(\cX)$ we write
\begin{equation}\begin{aligned} \label{eq:ip-rho}
 \ip{\Phi, \Psi}_{\rho} 
 & 
  := \frac12 \sum_{x,y \in \cX} 
         \Phi(x,y) \Psi(x,y) K(x,y) \rho(x,y) \pi(x)\;,
\\ \| \Phi \|_{\rho}
   & := \sqrt{\ip{\Phi, \Phi}_{\rho}}  \;,
\end{aligned}\end{equation}
and note that $\ip{\cdot, \cdot}_\pi = \ip{\cdot, \cdot}_\rho$ if $\rho(x) = 1$ for all $x\in \cX$.

For a probability density $\rho \in \cP(\cX)$ and $x \in \cX$ we consider the matrix $\hrho \in \R^{\cX \times\cX}$ defined by
\begin{align*}
  \hat\rho(x,y) := \rho(x,y)\;.
\end{align*}

Given two matrices $M, N \in \R^{\cX \times\cX}$, let $M \bullet N$ denote their entrywise product defined by
\begin{align*}
 (M \bullet N)(x,y) := M(x,y) N(x,y)
\end{align*}

The definition of $\cW$ can now be reformulated as follows:

\begin{lemma}[Geometric reformulation]\label{lem:metric-geom}
For $\bar\rho_0, \bar\rho_1 \in \cP(\cX)$ we have 
\begin{align*} 
\cW(\bar\rho_0, \bar\rho_1)^2
   = \inf_{\rho, \psi}  \bigg\{\int_0^1 \| \nabla\psi_t \|_{\rho_t}^2\dd t 
        \ : \ {(\rho, \psi) \in \CE(\bar\rho_0, \bar\rho_1)}  \bigg\}\;,
\end{align*}
and the differential equation in \eqref{eq:conditions} can be rewritten as
\begin{align}\label{eq:cont-eq}
 \dot\rho_t + \nabla \cdot (\hat\rho_t \bullet \nabla \psi_t) = 0\;.
\end{align}
\end{lemma}

\begin{proof}
This follows directly from the definitions.
\end{proof}

For the $L^2$-Wasserstein metric on Euclidean space, it is well known that one can take the infimum in the Benamou-Brenier formula \eqref{eq:Benamou-Brenier} over all vector fields $\Psi: \R^n \to \R^n$, rather than only considering gradients $\Psi = \nabla \psi$. 
In order to formulate a similar result in the discrete setting, we replace $(iv)$ and $(v)$ in \eqref{eq:conditions} by
 \begin{align} \label{eq:conditions-2} 
\begin{array}{ll}
{(iv')} & \Psi  : [0,T] \to \R^{\cX \times \cX} \text{ is measurable}\;;\\
{(v')} &  \text{For all $x \in \cX$ and a.e. $t \in (0,T)$ we have}\\
  &\displaystyle{\dot \rho_t(x) 
   + \frac12 \sum_{y \in \cX}
    \big(\Psi_t(x,y) - \Psi_t(y,x)\big) K(x,y) \rho_t(x,y) = 0}\;;\
\end{array}
\end{align}
and define
\begin{align*}
 \CE'(\rho_0, \rho_1) 
   := \{ (\rho, \Psi)  \ : \ (i), (ii), (iii), (iv'), (v') \text{ hold\,}
            \}\;.
\end{align*}
With this notation the following result holds.

\begin{lemma}\label{lem:only-grads}
For $\brho_0, \brho_1 \in \cP(\cX)$ we have 
\begin{align*}
 \cW(\bar\rho_0, \bar\rho_1)^2 
  = \inf \bigg\{ &  \frac12   \int_0^1 
  \sum_{x,y\in \cX} \Psi_t(x,y)^2
    		 K(x,y) \rho_t(x,y) \pi(x)
      \dd t 
         \ : \\ & \ {(\rho, \Psi) \in \CE'(\bar\rho_0, \bar\rho_1)}  \bigg\}\;.
\end{align*}
\end{lemma}

\begin{proof}
As the inequality ``$\geq$'' is trivial, it suffices to prove the inequality ``$\leq$''.
For this purpose, fix $\rho \in \cP(\cX)$ and let $\cH_\rho$ denote the set of all equivalence classes of functions $\Psi \in \R^{\cX \times \cX}$, where we identify functions that agree on $\{ (x,y) \in \cX \times \cX : \rho(x,y) K(x,y) > 0\}$. Endowed with the inner product $\ip{\cdot, \cdot}_\rho$ defined in \eqref{eq:ip-rho}, $\cH_\rho$ is a finite-dimensional Hilbert space. The discrete gradient $\nabla \phi(x,y) := \phi(x) - \phi(y)$  defines a linear operator $\nabla  : L^2(\cX,\pi) \to \cH_\rho$, whose adjoint is given by
\begin{align} \label{eq:adj}
 \nabla_\rho^* \Psi(x) := \frac12 \sum_{y \in \cX} \big(\Psi(x,y) - \Psi(y,x)  \big)
 K(x,y)   \rho(x,y)\;.
\end{align}
Let $P_\rho$ denote the orthogonal projection in $\cH_\rho$ onto the range of $\nabla$.

Now suppose that $((\rho_t),(\Psi_t)) \in \CE'(\brho_0, \brho_1)$ and let $\psi : [0,1] \to \R^\cX$ be such that $P_{\rho_t} \Psi_t = \nabla \psi_t$ for $t \in [0,1]$. In view of the orthogonal decomposition
\begin{align} \label{eq:orth-dec}
 \cH_\rho
   = \Ran(\nabla) \oplus^\perp \Ker(\nabla^*_{{\rho_t}})\;,
\end{align}
it follows that $(I - P_{\rho_t}) \Psi_t \in \Ker(\nabla^*_{{\rho_t}})$. This implies that $\nabla_{\rho_t}^* \Psi_t = \nabla_{\rho_t}^* (\nabla \psi_t)$, hence $(\rho, \psi) \in \CE(\brho_0, \brho_1)$.
Using the decomposition \eqref{eq:orth-dec} once more, we infer that 
 $\ip{\nabla \psi_t, \nabla \psi_t}_{\rho_t}
   \leq \ip{\Psi_t, \Psi_t}_{\rho_t}$, from which the result follows.
\end{proof}

\begin{remark}[Distance between positive measures]\label{rem:pos-dens}
It is of course possible, and occasionally useful, to extend the definition of  $\cW(\rho_0, \rho_1)$ to densities $\rho_0, \rho_1 : \cX \to \R_+$ having equal mass $m = \sum_{x \in \cX}\rho_i(x) \pi(x) \in (0,\infty) \setminus \{1\}$. A straightforward argument based on Lemma \ref{lem:only-grads} and the doubling property (A7) shows that 
\begin{align*}
c \cW( \rho_0, \rho_1)  \leq  \cW(\frac1m \rho_0, \frac{1}{m}\rho_1) \leq C \cW( \rho_0, \rho_1)\;,
\end{align*}
where the constants $c, C > 0$ do not depend on $\rho_0$ and $\rho_1$.
If (A7$'$) holds, it follows that
 $\cW( \rho_0, \rho_1) = {\sqrt{m}}  \cW(\frac1m \rho_0, \frac{1}{m}\rho_1)$.
\end{remark}

\subsection*{Basic properties of the metric}

The main result of this subsection reads as follows:

\begin{theorem}\label{thm:pseudo}
The mapping $\cW : \cP(\cX) \times \cP(\cX) \to \R$ defines a pseudo-metric on $\cP(\cX)$.
\end{theorem}

To prove this result we need some lemmas.

\begin{lemma} \label{lem:time-change}
For $\bar\rho_0, \bar\rho_1 \in \cP(\cX)$ and $T > 0$ we have
\begin{align*}
 \cW(\bar\rho_0, \bar\rho_1) 
&= \inf \bigg\{  \int_0^T [ A(\rho_t) \psi_t, \psi_t]^{\frac12}\dd t 
         \ : \ {(\rho, \psi) \in \CE_T(\bar\rho_0, \bar\rho_1)}  \bigg\}\;.
\end{align*}
\end{lemma}

\begin{proof}
This follows from a standard argument based on parametrisation by arc-length. We refer to \cite[Lemma 1.1.4]{AGS08} or \cite[Theorem 5.4]{DNS09} for the details in a very similar situation.
\end{proof}

The next lemma  provides a lower bound for $\cW$ in terms of the total variation distance, defined for $\rho_0, \rho_1 \in \cP(\cX)$ by
\begin{align*}
  d_{TV}(\rho_0,\rho_1) = \sum_{x \in \cX} \pi(x) |\rho_0(x) - \rho_1(x)|\;.
\end{align*}

\begin{lemma}[Lower bound by total variation distance]\label{lem:Lipschitz}
For $\rho_0, \rho_1 \in \cP(\cX)$ we have
\begin{align*}
 d_{TV}(\rho_0,\rho_1) \leq \sqrt{2 \|\theta \|_\infty} \cW(\rho_0, \rho_1)\;,
\end{align*}
where 
\begin{align*}
 \|\theta\|_{\infty} 
 = 
 \sup \big\{ \theta(s,t) : 0 \leq s,t \leq \big(\min_{x\in \cX} \pi(x)\big)^{-1}\big\}\;.
\end{align*}
\end{lemma}

\begin{proof}
We assume that $\cW(\rho_0, \rho_1) < \infty$, since otherwise there is nothing to prove.
Let $\eps > 0$, let $\rho_0, \rho_1 \in \cP(\cX)$ and take $(\rho,\psi) \in \CE(\rho_0,\rho_1)$ satisfying
\begin{align} \label{eq:almost}
  \int_0^1 [A(\rho_t)\psi_t, \psi_t] \dd t < \cW^2(\rho_0, \rho_1) + \eps\;.
\end{align}
Using the continuity equation \eqref{eq:cont-equ} we obtain for any $\phi : \cX \to \R$,
\begin{align*}
 \Big| &\sum_{x \in \cX} \phi(x) (\rho_0(x) - \rho_1(x)) \pi(x)\Big|
  = \bigg|\int_0^1 [\Pi \phi, \dot\rho_t] \dd t\bigg|
 \\ & =   \bigg|\int_0^1 [\Pi \phi, B(\rho_t) \psi_t]\dd t\bigg|
 = \bigg|\int_0^1 [A(\rho_t) \phi, \psi_t] \dd t\bigg|
 \\& 
   \leq 
   \bigg( \int_0^1 
    [A(\rho_t) \psi_t, \psi_t]  \dd t \bigg)^{1/2}
    \bigg( \int_0^1 
   [A(\rho_t) \phi, \phi] \dd t \bigg)^{1/2}\;,
\end{align*}
where the appeal to the Cauchy-Schwarz inequality is justified by \eqref{eq:A-nonneg}.
The latter integrand can be estimated brutally by
\begin{align*}
   [A(\rho_t) \phi, \phi]
	&= \frac12 \sum_{x,y \in \cX}
		 (\phi(x) - \phi(y))^2 K(x,y) \rho_t(x,y) \pi(x)
  \\&  \leq  
{2 \| \theta \|_\infty }
\| \phi \|_\infty^2 \sum_{x,y \in \cX} K(x,y) \pi(x)
 ={ 2 \| \theta \|_\infty}
\| \phi \|_\infty^2\;,
\end{align*}
where we used the stationarity of $\pi$ to obtain the latter identity.
Taking \eqref{eq:almost} into account, and noting that $\eps > 0$ is arbitrary, we thus obtain
\begin{align*}
 \Big| \sum_{x \in \cX} \phi(x) (\rho_0(x) - \rho_1(x)) \pi(x)\Big|
  \leq \sqrt{ 2 \| \theta \|_\infty } \|\phi\|_\infty  \cW(\rho_0, \rho_1)\;.
\end{align*}
Using the duality between $\ell^1(\cX)$ and $\ell^\infty(\cX)$, the result follows.
\end{proof}

\begin{proof}[Proof of Theorem \ref{thm:pseudo}]
The symmetry of $\cW$ is obvious, and Lemma \ref{lem:Lipschitz} implies that $\cW(\rho_0, \rho_1) > 0$ whenever $\rho_0 \neq \rho_1$. Finally, the triangle inequality easily follows using Lemma \ref{lem:time-change}.
\end{proof}

\subsection*{Characterisation of finiteness}

In the study of finiteness of the metric $\cW$, a crucial role will be played by the quantity
\begin{align*}
 C_\theta := \int_0^1 \frac{1}{\sqrt{\theta(1-r, 1+r)}} \dd  r \in [0,\infty]\;.
\end{align*}
Note that $C_\theta = \sqrt{2} \phi(1)$, where $\phi$ denotes the function defined in \eqref{eq:def-phi} with $p = q = 1$. Therefore $C_\theta$ is finite if and only if Dirac measures on the two-point space lie at finite $\cW$-distance from the uniform measure. Observe that $C_\theta < \infty$ if (A7$'$) holds, since in that case
\begin{align*}
 \theta(1-r,1+r)
  \geq \theta(1-r,1-r)
     = (1-r) \theta(1,1)\;,
\end{align*}
for $t \in [0,1)$.

The next result provides a characterisation of finiteness of the metric in terms of the support of the densities. For $\rho \in \cP(\cX)$ we shall write
\begin{align*}
 \supp \rho := \{ x\in \cX \ : \ \rho(x) > 0\}\;.
\end{align*}
Before stating the result we recall the following definition:

\begin{definition}\label{def:equiv}
Let $\rho \in \cP(\cX)$. For $x , y \in \cX$ we write `$x \sim_\rho y$'  if  
\begin{enumerate}[(i)]
\item $x = y$; or,
\item there exist $k \geq 1$ and $x_1, \ldots, x_k \in \cX$ such that
\begin{align*}
 \rho(x,x_1)  K(x,x_1), 
 \rho(x_1,x_2)K(x_1,x_2),
\;\ldots\;,
\rho(x_k,y)K(x_k,y) > 0\;.
\end{align*}

\end{enumerate}
\end{definition}

It is easy to see that for each $\rho \in \cP(\cX)$, $\sim_\rho$ defines an equivalence relation on $\cX$, which depends only on the support of $\rho$. Furthermore, 
if $\rho$ is strictly positive, then $x \sim_\rho y$ for any $x, y \in \cX$, since $K$ is irreducible by assumption.

Now we are ready to state the main result of this subsection.

\begin{theorem}[Characterisation of finiteness]\label{thm:finiteness}\mbox{}
\begin{enumerate}
\item 
If $C_\theta < \infty$, then $\cW(\rho_0, \rho_1) < \infty$ for all $\rho_0, \rho_1 \in \cP(\cX)$.
\item 
If $C_\theta = \infty$, the following assertions are equivalent for $\rho_0, \rho_1 \in \cP(\cX)$:
\begin{enumerate}
\item $\cW(\rho_0, \rho_1) < \infty$\;;
\item For any $x\in \cX$ we have
\begin{align} \label{eq:mass-conserv}
 \sum_{y\sim_{\rho_0} x} \rho_0(y) \pi(y) = 
 \sum_{y\sim_{\rho_1} x} \rho_1(y) \pi(y)\;.
\end{align}
\end{enumerate}
\end{enumerate}
\end{theorem}

Before turning to the proof of this result we record some immediate consequences:

\begin{corollary} \label{cor:finiteness}
Suppose that $C_\theta = \infty$. For $\rho_0, \rho_1 \in \cP(\cX)$ the following assertions hold:
\begin{enumerate}
\item If $\cW(\rho_0, \rho_1) < \infty$, then $\supp \rho_0 = \supp \rho_1$.
\item If $\supp \rho_0 = \supp \rho_1 = \cX$, then $\cW(\rho_0, \rho_1) < \infty$.
\end{enumerate}
\end{corollary}

\begin{proof}
(1) Suppose that $\rho_0(x) = 0$ for a certain $x \in \cX$. In view of (A5) it then follows that $x \not\sim_{\rho_0} y$ for any $y \neq x$, hence by Theorem \ref{thm:finiteness},
\begin{align*}
 \rho_1(x) \pi(x) \leq \sum_{y\sim_{\rho_1} x} \rho_1(y) \pi(y) =
 \sum_{y\sim_{\rho_0} x} \rho_0(y) \pi(y) =  \rho_0(x) \pi(x) = 0\;.
\end{align*}
It follows that $\rho_1(x) = 0$, which shows that $\supp \rho_0 \supseteq \supp \rho_1$. The reverse inclusion follows by reversing the roles of $\rho_0$ and $\rho_1$.

(2) If $\supp \rho_0 = \supp \rho_1 = \cX$, then $x \sim_{\rho_i} y$ for every $y \neq x$ and $i =0,1$ by irreducibility. It follows that
\begin{align*}
 \sum_{y\sim_{\rho_0} x} \rho_0(y) \pi(y) = 1 = 
 \sum_{y\sim_{\rho_1} x} \rho_1(y) \pi(y)\;,
\end{align*}
hence $\cW(\rho_0, \rho_1) < \infty$ by Theorem \ref{thm:finiteness}.
\end{proof}

The proof of Theorem \ref{thm:finiteness} relies on a sequence of lemmas of independent interest.

First we prove two comparison results, which relate the pseudo-metric $\cW$ on $\cP(\cX)$ to the pseudo-metric $\cW_{p,q}$ on $\cP(\cY)$, where $\cY = \{ a, b\}$ is a two-point space endowed with the Markov kernel \eqref{eq:K} with parameters $p$ and $q$.

\begin{lemma}[Comparison to the two-point space I]\label{lem:comparison-one}
Let $a, b \in \cX$ be distinct points with $K(a,b) >0$, and set $p := K(a,b) \pi(a)$.  Suppose that $\rho_0, \rho_1 \in \cP(\cX)$ satisfy $\rho_0(x) = \rho_1(x)$ for all $x \in \cX \setminus \{a, b\}$. 
Consider the two-point space $\cQ^1 = \{ \alpha, \beta\}$ endowed with the Markov kernel defined by $K(\alpha, \beta) := K(\beta, \alpha) := p$.
For $i = 0,1$, let $\brho_i: \cQ^1 \to \R_+$ be defined by
\begin{align*}
  \brho_i(\alpha) := 2 \rho_i(a) \pi(a)\;, \qquad
  \brho_i(\beta)  := 2 \rho_i(b) \pi(b)\;.
\end{align*}
Then we have
\begin{align*}
 \cW(\rho_0, \rho_1) \leq \sqrt{C_d}
 \cW_{p,p}(\brho_0, \brho_1)\;.
\end{align*}
where $C_d$ is the constant from (A7). In particular, if (A$7'$) holds, then
\begin{align*}
 \cW(\rho_0, \rho_1) \leq
 \cW_{p,p}(\brho_0, \brho_1)\;.
\end{align*}
\end{lemma}

\begin{remark}\label{rem:equal-mass}
Note that $\brho_0$ and $\brho_1$ are not necessarily probability densities on $\{\alpha, \beta\}$, but they do have equal mass, since
\begin{align*}
 \brho_i(\alpha)\pi(\alpha) +  \brho_i(\beta)\pi(\beta) 
   = \rho_j(a)\pi(a) + \rho_j(b)\pi(b)
\end{align*}
for $i, j \in \{0,1\}$. Therefore $\cW_{p,p}(\brho_0, \brho_1)$ can be interpreted in the sense of Remark \ref{rem:pos-dens}.
\end{remark}

\begin{proof}[Proof of Lemma \ref{lem:comparison-one}]
Let $\eps > 0$ and take $(\brho, \bar\psi) \in \CE( \brho_0, \brho_1)$.
It then follows that 
\begin{equation}\begin{aligned}\label{eq:ce2pt}
\dot \brho_t(\alpha) 
   +  ( \bpsi_t(\beta) - \bpsi_t(\alpha) )K(\alpha,\beta) \brho_t(\alpha,\beta) = 0\;,\\
\dot \brho_t(\beta) 
   +  ( \bpsi_t(\alpha) - \bpsi_t(\beta) )K(\beta,\alpha) \brho_t(\alpha,\beta) = 0\;.
\end{aligned}\end{equation}
For $t \in (0,1)$ define $\rho_t \in \cP(\cX)$ by
\begin{align*}
 \rho_t (a) := \frac{\brho_t(\alpha)}{2 \pi(a)}\;, \qquad
 \rho_t (b) :=  \frac{\brho_t(\beta)}{2 \pi(b)}\;, \qquad
 \rho_t (x) := \rho_0(x)\;,
\end{align*}
for $x \in \cX \setminus \{a,b\}$. 
Furthermore, we define $\Psi_t : \cX \times \cX \to \R$ by
\begin{align*}
 \Psi_t(a,b) &:= - \Psi_t(b,a) 
 			 := 
			    \frac{\brho_t(\alpha,\beta)}{2\rho_t(a,b)}
			    \big(\bpsi_t(\beta) - \bpsi_t(\alpha)\big)
			     \one_{\{\rho_t(a,b) > 0 \}}\;,\\
 \Psi_t(x,y) &:= 0\;,
\end{align*}
for all other values of $x,y \in \cX$.
Using \eqref{eq:ce2pt} it then follows that $(\rho, \Psi) \in \CE'(\rho_0, \rho_1)$. Using Lemma \ref{lem:only-grads} we thus obtain 
\begin{align*}
 \cW(\rho_0, \rho_1)^2 
   & \leq \int_0^1 \Psi_t(a,b)^2 \rho_t(a,b)  K(a,b) \pi(a)\dd t
 \\& = \frac12
      \int_0^1 \big(\bpsi_t(\alpha) - \bpsi_t(\beta)\big)^2
          \frac{\brho_t(\alpha,\beta)^2}{\rho_t(a,b)} 
          \one_{\{\rho_t(a,b) > 0 \}}
            K(\alpha,\beta) \pi(\alpha)\dd t\;.
\end{align*}
Using (A6) and (A7) we infer that
\begin{align*} 
 \brho_t(\alpha, \beta)
  & = \theta\big(2 \pi(a)\rho_t(a),2 \pi(b) \rho_t(b)\big)
 \leq 2C_d \theta\big( \rho_t(a),\rho_t(b)\big)
= 2 C_d\rho_t(a,b)\;,
\end{align*}
which yields
\begin{align*}
 \cW(\rho_0, \rho_1)^2 
   \leq C_d
    \int_0^1 \big(\bpsi_t(\alpha) - \bpsi_t(\beta)\big)^2
          \brho_t(\alpha,\beta)
            K(\alpha,\beta) \pi(\alpha)\dd t\;.
\end{align*}
Minimising the right-hand side over all $(\brho, \bar\psi) \in \CE( \brho_0, \brho_1)$, the result follows. 
\end{proof}

\begin{lemma}[Comparison to the two-point space II]\label{lem:comparison}
Let $\rho_0, \rho_1 \in \cP(\cX)$ and set $\beta_i(x) =  1- 2 \rho_i(x)\pi(x)$ for $i=0,1$ and $x \in \cX$. Then the bound
\begin{align*}
 \cW(\rho_0, \rho_1) \geq c \sup_{x \in \cX} \cW_{1,1}(\rho^{\beta_0(x)}, \rho^{\beta_1(x)})
\end{align*}
holds, for some $c > 0$ depending only on $K$, $\pi$ and $\theta$. 
If (A$7'$) and (A8) hold, then
\begin{align*}
 \cW(\rho_0, \rho_1) \geq \sup_{x \in \cX} \cW_{1,1}(\rho^{\beta_0(x)}, \rho^{\beta_1(x)})\;.
\end{align*}
\end{lemma}

\begin{proof}
First we shall  prove the result under the assumption that (A7$'$) and (A8) hold.
Fix $o \in \cX$ and let $\cY = \{a,b\}$ be a two-point space endowed with the Markov kernel \eqref{eq:K} with $p = q = 1$. For $\rho \in \cP(\cX)$ and $\psi \in \R^\cX$ we define, by a slight abuse of notation, $\rho \in \cP(\cY)$ and $\psi \in \R^\cY$ by
\begin{align*}
 \rho(a) &:= 2 \rho(o) \pi(o)\;,\qquad
 \rho(b) := 2 \sum_{x \neq o} \rho(x) \pi(x)\;, \\
 \psi(a) &:= \psi(o)\;,\qquad
 \psi(b) := \frac{\sum_{x \neq o} \psi(x)K(o,x)\rho(o,x)}
 		{\sum_{x \neq o}K(o,x)\rho(o,x)}\;.
\end{align*}
In the definition of $\psi(b)$ we use the convention that $0/0 = 0$.
Observe that $\rho$ indeed belongs to $\cP(\cY)$ since $\pi(a) = \pi(b) = \frac12$ and $\rho(a) + \rho(b) = 2$.
We set $\tilde\rho(a,b) := 2 \pi(o) \sum_{x \neq o} K(o,x)\rho(o,x)$ and claim that
\begin{align}
\label{eq:log-mean-bound}
 \tilde\rho(a,b)  
 &\leq \rho(a,b)\;,\\ 
\label{eq:cost-bound}
 [ A(\rho) \psi , \psi ] 
 &  \geq \frac12 (\psi(a) - \psi(b))^2 \tilde\rho(a,b)\;.
\end{align}
In the proof of both claims we shall assume that $\tilde\rho(a,b) > 0$, since otherwise there is nothing to prove.
To prove \eqref{eq:log-mean-bound}, note first that for any $x \in \cX$ with $K(o,x) > 0$, 
\begin{align} \label{eq:K-pi-ineq}
\frac{\pi(x)}{\pi(o)}
  = \frac{K(o,x)}{K(x,o)}
  \geq K(o,x)\;.
\end{align}
Using this inequality together with (A6), (A7') and (A8),
\begin{align*}
 \rho(a, b) 
  &  = \theta\bigg(   2 \rho(o) \pi(o),  2 \sum_{x \neq o} \rho(x) \pi(x)
    			\bigg)
  \\&  =  2 \pi(o)\theta \bigg(   \rho(o) ,   \sum_{x \neq o} \rho(x) \frac{\pi(x)}{\pi(o)}
    			\bigg)
  \\&  \geq  2 \pi(o)\theta \bigg(   \rho(o) ,   \sum_{x \neq o} K(o,x) \rho(x)      			\bigg)
  \\&  \geq  2 \pi(o) \sum_{x \neq o} K(o,x) \theta( \rho(o) , \rho(x)  )
       = \tilde\rho(a,b)\;,
\end{align*}
which proves \eqref{eq:log-mean-bound}.

To prove \eqref{eq:cost-bound}, write $k(x) := K(o,x)\rho(o,x)$ for brevity and note that
\begin{align*}
  \sum_{x\neq o } \psi(x)^2 k(x)
   \geq \frac{\big(\sum_{x\neq o } \psi(x) k(x)\big)^2}
            		  {\sum_{x\neq o } k(x)}
	  = 	\frac{\psi(b)^2 \tilde\rho(a,b)}{2 \pi(o)}\;.
\end{align*} 
Using the detailed balance equations \eqref{eq:detailed-balance} in the first inequality, we obtain
\begin{align*}
 [ A(\rho) \psi , \psi ] 
& =\frac12\sum_{x,y \in \cX} (\psi(x) - \psi(y))^2 K(x,y)\rho(x,y)\pi(x)
\\&  \geq \sum_{x\neq o } (\psi(o) - \psi(x))^2 K(o,x)\rho(o,x)\pi(o)
\\&  = \bigg( \psi(o)^2 \sum_{x\neq o }  k(x)
            - 2 \psi(o) \sum_{x\neq o } \psi(x) k(x)
            + \sum_{x\neq o } \psi(x)^2 k(x)\bigg)\pi(o)
\\& \geq \frac12\psi(a)^2 \tilde\rho(a,b)
           -  \psi(a)\psi(b)\tilde\rho(a,b)
            +  \frac12 \psi(b)^2 \tilde\rho(a,b) 
\\& = \frac12 (\psi(a) - \psi(b))^2 \tilde\rho(a,b)\;,
\end{align*}
which proves \eqref{eq:cost-bound}.

Take $(\rho, \psi) \in \CE(\rho_0, \rho_1)$. 
Since 
\begin{align*}
\dot\rho_t(o) + \sum_{x \neq o} (\psi_t(x) - \psi_t(o)) K(o,x) \rho_t(o,x) = 0\;,
\end{align*}
it follows that
\begin{align} \label{eq:cont2point-a}
  \dot\rho_t(a) + (\psi_t(b) - \psi_t(a)) \tilde \rho_t(a,b) = 0\;.
\end{align}
Set $\beta_t :=  1- 2 \rho_t(o) \pi(o)$ for $t \in [0,1]$ and note that $\dot\beta_t = 0$ if $\tilde \rho_t(a,b) = 0$.
Using \eqref{eq:cost-bound}, \eqref{eq:cont2point-a}, \eqref{eq:log-mean-bound} and Lemma \ref{lem:var-form} we obtain
\begin{align*}
  \int_0^1 [ A(\rho_t) \psi_t , \psi_t ] \dd t
 &  \geq \frac12 \int_0^1(\psi_t(a) - \psi_t(b))^2 \tilde\rho_t(a,b) \dd t
 \\& = \frac12 \int_0^1\frac{\dot\beta_t^2 \one_{\{ \tilde\rho_t(a,b) > 0 \}}}{\tilde\rho_t(a,b)} \dd t
     \geq \frac12 \int_0^1\frac{\dot\beta_t^2 \one_{\{ \rho_t(a,b) > 0 \}}}{\rho_t(a,b)} \dd t
\\&   \geq \cW_{1,1}^2(\rho^{\beta_0}, \rho^{\beta_1})\;.
\end{align*}
Taking the infimum over all pairs $(\rho, \psi) \in \CE(\rho_0, \rho_1)$, we infer that
\begin{align*}
 \cW^2(\rho_0, \rho_1) \geq \cW_{1,1}^2(\rho^{\beta_0}, \rho^{\beta_1})\;.
\end{align*}
The result follows by taking the supremum over $o \in \cX$.

Finally, without assuming (A7$'$) and (A8), the same argument applies, if one replaces \eqref{eq:log-mean-bound} by the following estimate, which uses 
the doubling property (A7), \eqref{eq:K-pi-ineq} and (A5):
\begin{align*}
\tilde\rho(a,b)
&  =   2 \pi(o) \sum_{x \neq o} K(o,x) \theta( \rho(o) , \rho(x)  )
  \\&  \leq C \sum_{x \neq o} 
  \theta\big(2\rho(o) K(o,x) \pi(o), 2 \rho(x) K(o,x) \pi(o)\big)
  \\& \leq C \sum_{x \neq o} 
  \theta\big(2\rho(o)  \pi(o), 2 \rho(x)  \pi(x)\big)
  \\& \leq C |\cX| 
  \theta\bigg(2\rho(o)  \pi(o), 2 \sum_{x \neq o} \rho(x)  \pi(x)\bigg)
  \\&  =   C |\cX| 
 \rho(a, b)\;.
\end{align*}
\end{proof}

The next lemma provides a useful characterision of the kernel and the range of the matrices $A(\rho)$ and $B(\rho)$.

\begin{lemma}\label{lem:Ran-Ker-B-1}
For $\rho \in \cP(\cX)$ we have
 \begin{align*}
\Ker A(\rho) = \Ker B(\rho) & = \{ \psi \in  \R^\cX \ | \ 
\psi(x) = \psi(y) \text{ whenever } x \sim_\rho y
		\}\;,\\
\Ran A(\rho) & =
 \Big\{ \psi  \in  \R^\cX  \ | \ \forall x \in \cX \ : \
 				    \sum_{y \sim_\rho x}  \psi(y) = 0	\Big\}\;,\\
\Ran B(\rho) & =
 \Big\{ \psi  \in  \R^\cX  \ | \ \forall x \in \cX \ : \
 				    \sum_{y \sim_\rho x}  \psi(y) \pi(y) = 0	\Big\}\;.
\end{align*}
\end{lemma}

\begin{proof}
Recall that (A3) and (A5) imply that $\rho(x,y) = 0$ whenever $\rho(x) = 0$ or $\rho(y) = 0$. Therefore the assertions concerning $A(\rho)$ follow directly from Lemma \ref{lem:matrix-kernel}.
Since $B(\rho) = \Pi^{-1} A(\rho)$, one has
\begin{align*}
 \Ker B(\rho) &= \Ker A(\rho)\;, \qquad \Ran B(\rho) = \Pi^{-1} \Ran A(\rho)\;,
\end{align*}
hence the remaining assertions follow as well.
\end{proof}

For $\sigma \in \cP(\cX)$ and $a \geq 0$ we shall use the notation
\begin{align*}
 \cP_\sigma^a(\cX) := \big\{ \rho \in \cP(\cX)  \ | \ &
 \forall x \in \cX \, : \, \eqref{eq:mass-conserv} \text{ holds with $\rho_0 = \rho$ and $\rho_1 = \sigma$; }
\\ &   \forall z \in \supp(\sigma) \, : \, \rho(z) \geq a \big\}\;.
\end{align*}

\begin{lemma}\label{lem:Ran-Ker-A-2}
For $\rho \in \cP(\cX)$, $B(\rho)$ restricts to an isomorphism from $\Ran A(\rho)$ onto $\Ran B(\rho)$.
Moreover, for $\sigma \in \cP(\cX)$ and $a > 0$ there exist constants $0 < c < C < \infty$ such that the bound
\begin{align} \label{eq:two-sided}
c \| \psi \| \leq  \| B(\rho) \psi \| \leq C \| \psi \|
\end{align}
holds for all $\rho \in \cP_\sigma^a(\cX)$ and all $\psi \in \Ran(\sigma)$.
\end{lemma}

\begin{proof}
Since $A(\rho)$ is self-adjoint, $A(\rho)$ restricts to an isomorphism  on its range. Since $\Pi$ is an isomorphism from $\Ran A(\rho)$ onto $\Ran B(\rho)$ and $B(\rho) = \Pi^{-1} A(\rho)$, the first assertion follows. 

Lemma \ref{lem:Ran-Ker-B-1} implies that $\Ran A(\rho) = \Ran A(\sigma)$ and $\Ran B(\rho) = \Ran B(\sigma)$ for all $\rho \in \cP_\sigma(\cX)$. Thus $B(\rho)$ restricts to an isomorphism, denoted by $B_\rho$, from $\Ran A(\sigma)$ onto $\Ran B(\sigma)$. Since the mapping $\cP_\sigma^a(\cX) \ni \rho \mapsto \| B_\rho^{-1}  \|$ is continuous w.r.t. the euclidean metric and strictly positive, the lower bound in \eqref{eq:two-sided} follows by compactness. The upper bound is clear, since the entries of $B(\rho)$ are bounded uniformly in $\rho$.
\end{proof}

The next result provides a partial converse to Lemma \ref{lem:Lipschitz}.

\begin{lemma}\label{lem:quant}
Fix $\sigma \in \cP(\cX)$ and $a > 0$.
There exist constants $0 < c < C < \infty$ such that for all $\rho_0, \rho_1 \in \cP_\sigma^a(\cX)$ 
we have 
\begin{align*}
c d_{TV}(\rho_0, \rho_1) \leq \cW(\rho_0, \rho_1) \leq C d_{TV}(\rho_0, \rho_1)\;.
\end{align*}
\end{lemma}

\begin{proof}
Since the lower bound for $\cW$ has been proved in Lemma \ref{lem:Lipschitz}, it remains to prove the upper bound.

For $t \in [0,1]$ set $\rho_t := (1-t) \rho_0 + t \rho_1$ and note that
$\rho_t \in \cP_\sigma^a(\cX)$.
Since 
\begin{align*}
\dot\rho_t = \rho_1 - \rho_0 \in  \Ran B(\rho_t) = \Ran B(\sigma)
\end{align*}
by Lemma \ref{lem:Ran-Ker-B-1}, Lemma \ref{lem:Ran-Ker-A-2} implies that, for each $t \in [0,1]$, there exists a unique element $\psi_t \in \Ran A(\rho_t)$ satisfying 
\begin{align*}
  \dot\rho_t = B(\rho_t)\psi_t\;.
\end{align*}
Moreover, Lemma \ref{lem:Ran-Ker-A-2} implies that
\begin{align*}
  \| \psi_t \| \leq C \| \rho_1 - \rho_0\|
\end{align*}
for some constant $C > 0$ that does not depend on $\rho_0, \rho_1$ and $t$. It thus follows that
\begin{align*}
 \cW(\rho_0, \rho_1)^2 \leq
  \int_0^1 [A(\rho_t) \psi_t, \psi_t] \dd t 
    \leq  C^2 C'\| \rho_1 - \rho_0 \|^2  
    \leq  C^2 C' C'' d_{TV}^2(\rho_0, \rho_1)\;,
\end{align*}
where $C' := \sup_{\rho \in \cP(\cX)} \|A(\rho)\| < \infty$ and $C'' > 0$ depends only on $\pi$.
\end{proof}

Now we are ready to prove the main result of this subsection.

\begin{proof}[Proof of Theorem \ref{thm:finiteness}]
Since $K$ is irreducible, (1) follows from Lemma \ref{lem:comparison-one}, Remark \ref{rem:pos-dens}  and the triangle inequality for $\cW$.

The implication $(b) \Rightarrow (a)$ of (2) follows from Lemma \ref{lem:quant}.

In order to prove the converse implication, we take $\rho_0, \rho_1 \in \cP(\cX)$ with $\cW(\rho_0, \rho_1) < \infty$ and claim that  $\supp \rho_0 = \supp \rho_1$. Indeed, if the claim were false, then there would exist $x \in \cX$ with $\rho_0(x) = 0$ and $\rho_1(x) > 0$ (or vice versa). Set $\beta = 1 - 2\pi(x)\rho_1(x)$ and note that $\beta \in [-1,1)$. Lemma \ref{lem:comparison} implies that $\cW(\rho_0, \rho_1) \geq c \cW_{1,1}(\rho^1, \rho^{\beta})$ for some $c > 0$. Since $C_\theta = \infty$, the right-hand side is infinite, which contradicts our assumption and thus proves the claim.

Let $(\rho, \psi) \in \CE(\rho_0, \rho_1)$ with $\int_0^1 [A(\rho_t), \psi_t, \psi_t] \dd t < \infty$. The claim implies that $\supp \rho_0 = \supp\rho_t$ for all $t \in [0,1]$ and therefore $x \sim_{\rho_t} y$ if and only if  $x \sim_{\rho_0} y$. Fix $z \in \supp \rho_0$ and take $x \in \cX$ with $x \sim_{\rho_0} z$. Since $K(x,y) \rho_t(x,y) = 0$ whenever $y \not\sim_{\rho_0} z$, we have
\begin{align*}
 \dot \rho_t(x) 
   + \sum_{y\sim_{\rho_t} z} ( \psi_t(y) - \psi_t(x) )K(x,y) \rho_t(x,y)
     				 = 0\;.
\end{align*}
Multiplying this identity by $\pi(x)$ and summing over $x \in \cX$ with $x \sim_{\rho_t} z$, it follows using the detailed balance equations \eqref{eq:detailed-balance} that 
\begin{align*}
 \sum_{x \sim_{\rho_0} z} \dot\rho_t(x) \pi(x) = 0\;,
\end{align*}
which implies \eqref{eq:mass-conserv}.
\end{proof}

\begin{remark}\label{rem:alternative}

Alternatively, the implication $(b) \Rightarrow (a)$ in the proof of Theorem \ref{thm:finiteness} can be proved as an application of Lemma \ref{lem:comparison-one}.

\end{remark}

We continue to prove the remaining parts of Theorem \ref{thm:main-metric}.

\begin{theorem}[Topology]\label{thm:topology}
Let $\sigma \in \cP(\cX)$. 
For $\rho, \rho_\alpha \in \cP_\sigma(\cX)$, the following assertions are equivalent:
\begin{align*}
 (1) \quad \lim_\alpha d_{TV}(\rho_ \alpha, \rho) = 0\;; \qquad
  (2) \quad\lim_\alpha \cW(\rho_ \alpha, \rho) = 0\;.
\end{align*}
\end{theorem}

\begin{proof}
It follows from Lemma \ref{lem:Lipschitz} that (2) implies (1). 

Conversely, suppose that (1) holds. If $C_\theta < \infty$, then (2) follows easily using Lemma \ref{lem:comparison-one}. If $C_\theta = \infty$, there exists an index $\bar \alpha$ and a constant $b > 0$ such that $\rho$ and $\rho_\alpha$ belong to $\cP_\sigma^b(\cX)$ for every $\alpha \geq \bar \alpha $. Lemma \ref{lem:quant} implies then that there exists a constant $C > 0$ such that 
\begin{align*}
 \cW(\rho_\alpha, \rho) \leq C d_{TV}(\rho_\alpha, \rho)
\end{align*}
for all $\alpha \geq \bar \alpha$, which yields the result.
\end{proof}

\begin{theorem}[Completeness]\label{thm:completeness}
For every $\sigma \in \cP(\cX)$ the metric space $(\cP_\sigma(\cX), \cW)$ is complete.
\end{theorem}

\begin{proof}
If $C_\theta < \infty$, this follows directly from Lemma \ref{lem:Lipschitz} and Theorem \ref{thm:topology}.
If $C_\theta = \infty$, take a sequence $(\rho_n)_n$ in $\cP_\sigma(\cX)$ which is Cauchy with respect to $\cW$. In particular,  $(\rho_n)_n$  is bounded in the $\cW$-metric, hence by Lemma \ref{lem:comparison} there exists a constant $a > 0$ such that $\rho_n$ belongs to $\cP_\sigma^a(\cX)$ for every $n$. By Lemma \ref{lem:Lipschitz} $(\rho_n)_n$ is Cauchy in the total variation metric, hence $\rho_n$ converges to some $\bar\rho \in \cP(\cX)$ in total variation. Since $\cP_\sigma^a(\cX)$ is a $d_{TV}$-closed subset of $\cP(\cX)$, it follows that $\bar\rho$ belongs to $\cP_\sigma^a(\cX)$.
From Theorem \ref{thm:topology} we then infer that $\rho_n$ converges to $\bar\rho$ in $\cW$-metric, which yields the desired result.
\end{proof}

\subsection*{Riemannian structure}

Fix a probability density $\sigma \in \cP(\cX)$ and consider the space

\begin{align*}
 \cP_\sigma'(\cX) := 
  \Big\{ \rho \in \cP(\cX) \ \Big| \ \forall x \in \cX \ : \ 
  \sum_{y\sim_{\rho} x} \rho(y) \pi(y) = 
 \sum_{y\sim_{\sigma} x} \sigma(y) \pi(y)
 \Big\}\;.
\end{align*}
Note that $\cP_{\one}'(\cX) = \cP_*(\cX)$ where $\one$ denotes the uniform density with respect to $\pi$.
Moreover, if $C_\theta = \infty$, Theorem \ref{thm:finiteness} implies that  $\cP_\sigma'(\cX) =  \cP_\sigma(\cX)$ for all $\sigma \in \cP(\cX)$.

Our next aim is to show that the metric space $(\cP_\sigma(\cX), \cW)$ is a   Riemannian manifold.
First, we have the following result:

\begin{proposition}\label{prop:smooth-mfd}
The metric space $(\cP_\sigma'(\cX), \cW)$ is a smooth manifold of dimension   
\begin{align*}
d(\sigma) := |\supp \sigma| - n(\sigma)\;,
\end{align*}
where $|\supp \sigma|$ is the cardinality of $\supp \sigma$, and $n(\sigma)$ is the number of equivalences classes in the support of $\sigma$ for the equivalence relation $\sim_\sigma$. 
\end{proposition}

\begin{proof}
It follows from Theorem \ref{thm:finiteness} and Lemma \ref{lem:Ran-Ker-B-1} that $\cP_\sigma'(\cX)$ is a relatively open subset of the affine subspace
\begin{align*}
  S_\sigma := \sigma + \Ran B(\sigma) \subseteq \R^\cX\;.
\end{align*}
Theorem \ref{thm:topology} implies that the topology induced by $\cW$ coincides with the euclidean topology on $\cP_\sigma'(\cX)$, hence $(\cP_\sigma'(\cX), \cW)$ is a smooth manifold.

The assertion concerning the dimension follows immediately, since $d(\sigma)$ is the dimension of $\Ran B(\sigma)$.
\end{proof}

Fix $\sigma \in \cP(\cX)$ and $\rho \in \cP_\sigma'(\cX)$. Since $\cP_\sigma'(\cX)$ is an open subset of the affine space $\sigma + \Ran B(\sigma)$, the tangent space of $\cP_\sigma'(\cX)$ at $\rho$ can be naturally identified with $\Ran B(\sigma) = \Ran B(\rho)$.
Our next aim is to show that the tangent space can be identified with a space of gradients, in the spirit of the Otto calculus developed in \cite{O01}. In fact, we shall construct an isomorphism $\cI_\rho$ from $\Ran B(\sigma)$ onto
\begin{align*}
T_\rho := \{ \nabla \psi \in \R^{\cX \times \cX} \ : \ 
               \psi \in \Ran A(\rho) \}\;.
\end{align*}

\begin{remark}\label{rem:psi-s}
Note that if $\rho$ belongs to $\cP_*(\cX)$, we have
\begin{align*}
 T_\rho = \{ \nabla \psi \in \R^{\cX \times \cX} \ : \ 
               \psi \in \R^\cX \}\;.
\end{align*}
However, it is easy to see that this is no longer true if $\rho \notin \cP_*(\cX)$.
\end{remark}

\begin{proposition}\label{prop:isom}
Let $\rho \in \cP_\sigma'(\cX)$. The mapping
\begin{align*}
 \cI_\rho : \Ran B(\sigma) \to T_\rho\;, \qquad
 B(\rho) \psi \mapsto \nabla \psi
\end{align*}
defined for $\psi \in \Ran A (\rho)$, is a linear isomorphism.
\end{proposition}

\begin{proof}
To show that $\cI_\rho $ is well-defined, consider the following mappings:
\begin{align*}
  F_\rho : \Ran A(\rho) &\to \Ran B(\rho)\;,& 
  \psi &\mapsto B(\rho) \psi\;,\\
  G : \Ran A(\rho) &\to T_\rho\;,& 
   \psi  &\mapsto \nabla\psi\;.
\end{align*}
We claim that $F_\rho$ and $G$ are linear isomorphisms. Once this has been established, the proposition follows at once.
The claim for $F_\rho$ has been proved in Lemma \ref{lem:Ran-Ker-A-2}. 
To prove the claim for $G$, suppose that $\nabla \psi = 0$ for some $\psi \in \Ran(A)$. It then follows that
\begin{align*}
 [A(\rho) \psi, \psi] = \ip{\nabla\psi, \nabla\psi}_\rho =0 , 
\end{align*}
Since $A(\rho)$ is symmetric and $\psi \in \Ran A(\rho)$, it follows that $\psi = 0$, which completes the proof.
\end{proof}

The following statement clarifies the connection with the Otto calculus in the continuous setting:

\begin{proposition}\label{prop:cont-eq}
Let $\rho : [0,1] \to \cP_\sigma'(\cX)$ be differentiable at $t \in [0,1]$. Then $\cI_{\rho_t} \dot \rho_t$ is the unique element $\nabla \psi_t \in T_{\rho_t}$ satisfying the identity
\begin{align*}
 \dot \rho_t + \nabla\cdot ( \hrho_t \bullet \nabla \psi_t) = 0\;.
\end{align*}
\end{proposition}

\begin{proof}
Since $B(\rho)\psi = - \nabla \cdot (\hrho \bullet\nabla \psi)$ for $\rho \in \cP(\cX)$ and $\psi \in \R^\cX$, this is an immediate consequence of Proposition \ref{prop:isom}.
\end{proof}

Henceforth we shall identify the tangent space of $\cP_\sigma'(\cX)$ at $\rho$ with $T_\rho$ by means of the isomorphism $\cI_\rho $.

\begin{definition}\label{def:ip}
Let $\rho \in \cP_\sigma'(\cX)$. We endow $T_\rho$ with the inner product
\begin{align*}
 \ip{\nabla \phi, \nabla \psi}_\rho 
   = \frac12 \sum_{x,y\in \cX} 
 (\phi(x) - \phi(y))(\psi(x) - \psi(y)) K(x,y) \rho(x,y)\pi(x)\;,
\end{align*}
defined for $\phi, \psi \in \Ran A(\rho)$.
\end{definition}

Note that, for $\rho \in \cP_\sigma'(\cX)$ and $\phi, \psi \in \Ran A(\rho)$
\begin{align} \label{eq:inp}
 \ip{\nabla \phi, \nabla \psi}_\rho 
 =  [A(\rho) \phi, \psi]\;.
\end{align}

\begin{remark}\label{rem:ip}
It is clear from the definition that $\ip{\nabla \phi, \nabla \psi}_\rho$ is well-defined. Moreover, \eqref{eq:inp} implies that if $\ip{\nabla \psi, \nabla \psi}_\rho = 0$ for some $\psi \in \Ran A(\rho)$, then $\psi = 0$, thus the expression indeed defines an inner product on $T_\rho$.
\end{remark}

\begin{theorem}\label{thm:RiemMfd} The following statements hold:
\begin{itemize}
\item If $C_\theta < \infty$ and (A8) holds, then $(\cP_*(\cX),\cW)$ is a Riemannian manifold.
\item If $C_\theta = \infty$, then $(\cP_\sigma'(\cX), \cW)$ is a complete Riemannian manifold for every $\sigma \in \cP(\cX)$.
\end{itemize}
The Riemannian metric is given by Definition \ref{def:ip}.

\end{theorem}

\begin{proof}
Suppose first that $C_\theta = \infty$. Then Proposition \ref{prop:smooth-mfd} asserts that $(\cP_\sigma(\cX), \cW)$ is a smooth manifold and the completeness has been proved in Theorem \ref{thm:completeness}.
The result would follow immediately from Lemma \ref{lem:metric-geom} and Definition \ref{def:ip}, if we were allowed to add the following requirements to the definition of $\CE(\rho_0, \rho_1)$ without changing the value of $\cW(\rho_0, \rho_1)$:
\begin{enumerate}[(i)]
\item $\rho_t \in \cP_\sigma(\cX)$ for all $t \in [0,1]$;
\item $\psi_t \in \Ran A(\rho_t)$ for all $t \in [0,1]$.
\end{enumerate}
But (i) may be added by Theorem \ref{thm:finiteness} and (ii) may be added in view of the orthogonal decomposition $\cX = \Ran A(\rho) \oplus \Ker A(\rho)$.

If $C_\theta < \infty$ the same argument applies, with Lemma \ref{lem:strictly-pos} below providing the analogue of (i).
\end{proof}

The next result asserts that in the definition of $\cW$, only curves consisting of strictly positive densities need to be considered if the endpoints are strictly positive as well.

\begin{lemma}\label{lem:strictly-pos}
Suppose that (A8) holds. For $\rho_0, \rho_1 \in \cP_*(\cX)$, we may replace $(iii)$ in Definition \ref{def:metric} by ``$(iii') : \rho_t \in \cP_*(\cX) \text{ for all $t \in [0,T]$''}.$
\end{lemma}

\begin{proof}
For notational reasons, let us write
\begin{align*}
  \cA(\rho, \Psi) := \| \Psi\|_\rho^2
   =  \frac12   
  \sum_{x,y\in \cX} \Psi(x,y)^2
    		 K(x,y) \rho(x,y) \pi(x)
\end{align*}
for $\rho \in \cP(\cX)$ and $\Psi \in \R^{\cX \times \cX}$.
Let $0 < \eps < 1$ and let $(\rho, \Psi) \in \CE'(\rho_0, \rho_1)$ be such that
\begin{align*}
\int_0^1   \cA(\rho_t, \Psi_t) 
      \dd t < \cW^2(\rho_0, \rho_1) + \eps\;.
\end{align*}
We set $\rho_i^\eps = (1 - \eps) \rho_i + \eps$ for $i =0,1$.

Firstly, we define $(\rho^\eps, \Psi^\eps) \in \CE'(\rho_0^\eps, \rho_i^\eps)$ by 
\begin{align*}
\rho_t^\eps(x) &:= (1-\eps)\rho_t(x) + \eps\;,\\
\Psi_t^\eps(x,y) &:= (1-\eps) 
     \frac{\rho_t(x,y)}{\rho_t^\eps(x,y)} 
 \Psi_t(x,y)\;.
\end{align*}
The concavity assumption (A8) implies the convexity of the  function 
\begin{align*}
  \R \times \R_+ \times \R_+ \ni 
     (x, s,t) \mapsto \frac{x^2}{\theta(s,t)}\;,
\end{align*}
which yields
\begin{align*}
 \int_0^1 \cA(\rho_t^\eps, \Psi_t^\eps)
      \dd t 
      \leq (1- \eps ) \int_0^1 \cA(\rho_t, \Psi_t) \dd t
      < (1-\eps) \cW^2(\rho_0, \rho_1) +  \eps\;.
\end{align*}

Secondly, for $i = 0,1$, we define $(\rho^{i,\eps}, \Psi^{i,\eps}) \in \CE'(\rho_i, \rho_i^\eps)$ by linear interpolation, i.e., 
\begin{align*}
 \rho_t^{i,\eps} := (1-t) \rho_i + t \rho_i^{\eps}\;.
\end{align*}
As in the proof of Lemma \ref{lem:quant}, for $t \in (0,1)$, let 
$\psi_t^{i,\eps}$ be the unique element in $\Ran A(\rho_t^{i,\eps})$ satisfying 
$\dot\rho_t^{i,\eps} = B(\rho_t^{i,\eps}) \psi_t^{i,\eps}$. Setting $\Psi^{i,\eps} := \nabla\psi^{i,\eps}$, it then follows that $(\rho^{i,\eps}, \Psi^{i,\eps}) \in \CE'(\rho_i, \rho_i^\eps)$.
Lemma \ref{lem:quant} and its proof imply that there exists a constant $C > 0$, independent of $\eps > 0$, such that
\begin{align*}
  \int_0^1 \cA(\rho_t^{i,\eps},  \Psi_t^{i,\eps}) \dd t 
    \leq C d_{TV}^2(\rho_i, \rho_i^\eps)
    \leq 4C \eps^2\;.
\end{align*}

Finally, it remains to rescale the three curves in time and glue them together. We thus define
\begin{equation*}
(\brho_t^\eps, \bPsi_t^\eps) := \begin{cases}
    \big(\rho_{t/\eps}^{0,\eps}, \eps^{-1} \Psi_{t/\eps}^{0,\eps}\big)\;,
    					    			& t \in [0, \eps]\;,\\
       \big(\rho_{(t - \eps)/(1-2\eps)}^\eps, 
       		(1-2\eps)^{-1}\Psi_{(t - \eps)/(1-2\eps)}^\eps\big)\;,
       							& t \in (\eps, 1 - \eps)\;,\\
	\big(\rho_{(1-t)/\eps}^{1,\eps}, \eps^{-1}\Psi_{(1-t)/\eps}^{1,\eps}\big)\;,
						        & t \in [1-\eps, 1]\;,
\end{cases}
\end{equation*}
so that $(\brho^\eps, \bPsi^\eps) \in \CE(\rho_0, \rho_1)$. We infer that
\begin{align*}
  \int_0^1 \cA(\brho_t^\eps, \bPsi_t^\eps) \dd t 
&    \leq   \int_0^1
    		    \frac{\cA(\rho_t^{0,\eps},  \Psi_t^{0,\eps})}{\eps}
         +  \frac{\cA(\rho_t^\eps, \Psi_t^\eps)}{1 - 2\eps}
         +  \frac{\cA(\rho_t^{1,\eps},  \Psi_t^{1,\eps})}{\eps} \; dt
\\&    \leq
    		        4C \eps 
         +  \frac{ (1-\eps) \cW^2(\rho_0, \rho_1) +  \eps}{1 - 2\eps}
         +  		4C \eps\;.
\end{align*}
Since the right-hand side tends to $\cW^2(\rho_0, \rho_1)$ as $\eps \to 0$, the result follows from the observation that $\bPsi_t^\eps$ may be replaced by $P_{\brho_t^\eps} \bPsi_t^\eps$, as in the proof of Lemma \ref{lem:only-grads}.
\end{proof}

In the next result we will slightly abuse notation and write
\begin{align*}
 \partial_1\rho(x,y)  := \partial_1 \theta( \rho(x), \rho(y))\;.
\end{align*}

\begin{theorem}[Geodesics]\label{thm:geodesics}

Suppose that $C_\theta = \infty$ and let $\sigma \in \cP(\cX)$. 
The following assertions hold:
\begin{enumerate}
\item For each $\bar\rho_0, \bar\rho_1 \in \cP_\sigma(\cX)$ there exists a constant speed geodesic $\rho : [0,1] \to \cP(\cX)$ with $\rho_0 = \bar\rho_0$ and $\rho_1 = \bar\rho_1$.
\item Let $\rho : [0,1] \to \cP_\sigma(\cX)$ be a constant speed geodesic and let $\psi_t = \cI_{\rho_t}\dot\rho_t$. Then the following equations hold for $t \in [0,1]$ and $x \in \cX$:
\begin{equation}\begin{aligned} \label{eq:geod-equs}
\begin{cases}
\partial_t \rho_t(x)  = 
   \displaystyle\sum_{y \in \cX}  ( \psi_t(x) - \psi_t(y) ) K(x,y) \rho_t(x,y)
     		 	 \;,\\
 \partial_t \psi_t(x) = \displaystyle\frac12
     \displaystyle\sum_{y \in \cX} \big(  \psi_t(x) -\psi_t(y) \big)^2 
     		 K(x,y) \partial_1\rho_t(x,y) 
   \;.
\end{cases}
\end{aligned}\end{equation}
\end{enumerate}
\end{theorem}  

\begin{proof}
Since $(\cP_\sigma(\cX), \cW)$ is a complete Riemannian manifold, (1) follows from the Hopf-Rinow theorem.
The equations in (2) are the equations for the cogeodesic flow (see, e.g., \cite[Theorem 1.9.3]{Jo08}) and follow directly from the representation of $\cW$ as a Riemannian metric given in this section.
\end{proof}

\begin{remark}\label{rem:compare}
The equations \eqref{eq:geod-equs} should be compared to the geodesic equations for the $L^2$-Wasserstein metric over $\R^n$ (see \cite{BB00}, \cite{O01}, \cite{OV00}), which are given under appropriate assumptions by
\begin{equation}\begin{aligned}\label{eq:geod-Rn}
 \left\{ \begin{array}{l}
 \partial_t \rho + \nabla \cdot (\rho \nabla \psi) = 0\;,\\
 \partial_t \psi + \frac12|\nabla \psi|^2  = 0 \;.\end{array} \right.
\end{aligned}\end{equation}
The equations \eqref{eq:geod-equs} are a natural discrete analogue of \eqref{eq:geod-Rn}. Note however that the equations for $\psi$ in the discrete case depend on $\rho$. 
\end{remark}

\section{Gradient flows of entropy functionals} 
  \label{sec:gradFlow}

We continue in the setting of Section \ref{sec:metric}, where $K$ is an irreducible and reversible Markov kernel on a finite set $\cX$. 
We fix a function $\theta : \R_+ \times \R_+ \to \R_+$ satisfying Assumption \ref{ass:theta-2} and consider the associated (pseudo-)metric defined in Section \ref{sec:metric}. If $C_\theta < \infty$, we shall also assume that (A8) holds.

Since $\cP_*(\cX)$ is a Riemannian manifold, as has been shown in Theorem  \ref{thm:RiemMfd}, we are in a position to study gradient flows of smooth functionals defined on $\cP_*(\cX)$. Let 
\begin{align*}
\Delta := K - I
\end{align*}
denote the generator of the continuous time Markov semigroup $(e^{t\Delta})_{t \geq 0}$ associated with $K$. 
The main result in this section is Theorem \ref{thm:grad-flow-Riem}, which asserts that solutions to the ``heat equation'' $\dot\rho_t = \Delta\rho_t$ are gradient flow trajectories of the entropy $\cH$ with respect to the metric $\cW$.

\begin{notation*}\label{not:various}
In view of Proposition \ref{prop:isom}, we shall always regard $T_\rho$ as being the tangent space of $\cP_*(\cX)$ at $\rho \in \cP_*(\cX)$. 
The tangent vector field along a smooth curve $t \mapsto \rho_t \in \cP_*(\cX)$ will be denoted by
\begin{align*}
 t \mapsto D_t \rho \in T_{\rho_t}\;.
\end{align*}
The gradient of a smooth functional $\cG : \cP_*(\cX) \to \R$ at $\rho \in \cP_*(\cX)$ is denoted by 
\begin{align*}
 \grad \cG(\rho) \in T_{\rho}\;.
\end{align*}
\end{notation*}

\subsection*{Functionals}
We shall consider the following types of functionals:
\begin{itemize}
\item
For a function $V : \cX \to \R$ we consider the \emph{potential energy functional} $\cV : \cP_*(\cX) \to \R$ defined by
\begin{align*}
 \cV(\rho) := \sum_{x \in \cX} V(x) \rho(x) \pi(x)\;.
\end{align*}
\item
For a differentiable function $f : (0,\infty) \to \R$, we consider the \emph{generalised entropy} 
$\cF : \cP_*(\cX) \to \R$ defined by 
\begin{align*}
 \cF(\rho) &:= \sum_{x \in \cX} f(\rho(x)) \pi(x)\;.
\end{align*} 
\end{itemize}

\begin{proposition}[Gradient of potential energy functionals]\label{prop:grad-V}
The functional $\cV : \cP_*(\cX) \to \R$ is differentiable, and for $\rho \in \cP_*(\cX)$ we have
\begin{align*}
 \grad \cV(\rho) = \nabla V \;.
\end{align*}
\end{proposition}

\begin{proof}
Clearly, $\cV$ is differentiable.
Let $t \mapsto \rho_t \in \cP_*(\cX)$ be a differentiable curve and let $\psi_t \in \Ran A(\rho_t)$ be such that $\nabla\psi_t := D_t \rho$. Then
\begin{align*}
\ddt\cV(u_t)
   &	 = \sum_{x\in \cX} V(x) \dot \rho_t(x)  \pi(x)
     =  \sum_{x\in \cX} V(x)
 			 (B(\rho_t) \psi_t)(x) \pi(x)
 \\& = - \ip{V, \nabla\cdot(\hrho_t \bullet\nabla \psi_t)}_\pi
     =   \ip{\nabla V, \hrho_t \bullet\nabla \psi_t}_\pi
     =   \ip{\nabla V, \nabla\psi_t}_{\rho_t}\;,
\end{align*}
which yields the result.
\end{proof}

\begin{proposition}[Gradient of generalised entropy functionals]\label{prop:grad-F}
The functional $\cF : \cP_*(\cX) \to \R$ is differentiable, and for $\rho \in \cP_*(\cX)$ we have
\begin{align*}
 \grad \cF(\rho) = \nabla (f' \circ \rho) \;.
\end{align*}
\end{proposition}

\begin{proof}
The differentiability of $\cF$ is clear from its definition.
Let $t \mapsto \rho_t \in \cP_*(\cX)$ be a differentiable curve and let $\psi_t \in \Ran A(\rho_t)$ be such that $\nabla\psi_t := D_t \rho$. Since $f$ is differentiable, we obtain
\begin{align*}
\ddt\cF(u_t)
   &	 = \sum_{x\in \cX} f'(\rho_t(x)) \dot \rho_t(x)  \pi(x)
     =  \sum_{x\in \cX} f'(\rho_t(x))
 			 (B(\rho_t) \psi_t)(x) \pi(x)
 \\& = - \ip{f'(\rho_t), \nabla\cdot(\hrho_t \bullet\nabla \psi_t)}_\pi
     =   \ip{\nabla f'(\rho_t), \hrho_t \bullet\nabla \psi_t}_\pi
 \\&    =   \ip{\nabla f'(\rho_t), \nabla\psi_t}_{\rho_t}\;,
\end{align*}
which yields the result.
\end{proof}

In the special case where $\cF = \cH$ is the entropy functional from \eqref{eq:entropy} we obtain:

\begin{corollary}\label{cor:grad-H}
The functional $\cH : \cP_*(\cX) \to \R$ is differentiable, and for $\rho \in \cP_*(\cX)$ we have
\begin{align*}
 \grad \cH(\rho) = \nabla \log \rho\;.
\end{align*}
\end{corollary}

\begin{proof}
This follows directly from Proposition \ref{prop:grad-F}. 
\end{proof}

\subsection*{Gradient flows}

In order to study gradient flows, we impose the following assumption which will be in force throughout the remainder of this section.

\begin{assumption}\label{ass:}
In addition to Assumption \ref{ass:theta-2} we assume:
\begin{itemize}
\item[(A9)] 
There exists a function $k \in C^1((0,\infty); \R)$
such that
\begin{align*}
 \theta(s,t) = \frac{s-t}{k(s) - k(t)}
\end{align*}
for all $s, t > 0$ with $s \neq t$.
\end{itemize}
\end{assumption}

Recall that this assumption is satisfied if $\theta$ is the logarithmic mean, in which case $k(t) = \log(t)$.

\begin{proposition}[Tangent vector field along the heat flow]\label{prop:ent-heat}
Let $\rho \in \cP(\cX)$ and let $\rho_t = e^{t\Delta} \rho$, $t \geq 0$ denote the heat flow. Then $t \mapsto \rho_t$ is $C^\infty$ on $(0,\infty)$ and
for $t > 0$ we have
\begin{align*}
 D_t \rho = - \nabla (k \circ \rho_t)\;.
\end{align*}
\end{proposition} 

\begin{proof}
The differentiability assertion follows from general Markov chain theory. 
For any $\rho \in \cP_*(\cX)$, we have
\begin{align*}
 \rho(x,y) = \frac{\rho(x) - \rho(y)}{k(\rho(x)) - k(\rho(y))}\;,
\end{align*}
and therefore
\begin{align*}
 \Delta \rho 
 = \nabla \cdot (\nabla \rho)
 = \nabla \cdot (\hrho \bullet \nabla (k \circ \rho) )\;.
\end{align*}
Since $t \mapsto\rho_t$ solves the heat equation $\dot\rho_t = \Delta \rho_t$,
it follows that 
\begin{align*}
\dot\rho_t - \nabla \cdot (\hrho_t \bullet \nabla ( k \circ \rho_t)) = 0\;,
\end{align*}
hence $D_t \rho = - \nabla (k \circ \rho_t)$ by Proposition \ref{prop:cont-eq}.
\end{proof}

We slightly modify the usual definition of a gradient flow trajectory, as we wish to allow for initial values that do not belong to $\cP_*(\cX)$:
 
\begin{definition}[Gradient flow]\label{def:grad-flow}
Let $\cF : \cP_*(\cX) \to \R$ be differentiable. A curve $\rho : [0,\infty) \to \cP(\cX)$ is said to be a \emph{gradient flow trajectory} for $\cF$ starting from $\bar\rho \in \cP(\cX)$ if the following assertions hold:
\begin{enumerate}
\item $t \mapsto \rho_t$ is differentiable on $(0,\infty)$, for every $t > 0$ we have $\rho_t \in \cP_*(\cX)$ and 
\begin{align*}
  D_t \rho = - \grad \cF(\rho_t)\;.
\end{align*}
\item $t \mapsto \rho_t$ is continuous in total variation at $t = 0$ and $\rho_0 = \bar\rho$.
\end{enumerate}
\end{definition}

\begin{theorem}\label{thm:grad-flow-Riem}
Let $f \in C^2((0,\infty) ; \R)$ be such that $f' = k$ and let $\rho \in \cP(\cX)$. Then the heat flow $t \mapsto e^{t\Delta} \rho$ is a gradient flow trajectory for the functional $\cF$ with respect to $\cW$.
\end{theorem}

\begin{proof}
The first condition in Definition \ref{def:grad-flow} is a consequence of  Propositions \ref{prop:grad-F} and \ref{prop:ent-heat}. The second one follows from general Markov chain theory.
\end{proof}

\begin{corollary}[Heat flow is gradient flow of the entropy]\label{cor:grad-flow-Riem}
Let $\theta$ be the logarithmic mean defined by 
$\theta(s,t) = \int_0^1 s^{1-p} t^p \dd p
$ and let $\rho \in \cP(\cX)$. Then the heat flow $t \mapsto e^{t\Delta} \rho$ is a gradient flow trajectory for the entropy $\cH$ with respect to $\cW$.
\end{corollary}

\begin{proof}
This is a special case of Theorem \ref{thm:grad-flow-Riem} with $k(t) = 1 + \log t$ and $f(t) = t \log t$.
\end{proof}

 \appendix

\section{A result from the theory of diagonally dominant matrices}

The following result from the theory of diagonally theory is a special case of \cite{Dahl00}. For the convenience of the reader we present a simple proof.

\begin{lemma}\label{lem:matrix-kernel}
Let $A = (a_{ij})_{i,j = 1, \ldots, n}$ be a real matrix satisfying
\begin{align*}
(1) \  \forall i \ : \ a_{ii} \geq 0\;,  \quad \ \
(2) \ \forall i \neq j  \ : \ a_{ij} = a_{ji} \leq 0\;,  \quad \ \
(3)\ \forall i \ : \ \sum_{j} a_{ij} = 0\;.
\end{align*}
Consider the equivalence relation $\sim$ on $I = \{1, \ldots, n\}$ defined by
\begin{align*}
 i \sim j \quad : \Leftrightarrow \quad
  \left\{ \begin{array}{l}
	i = j\;, \qquad \text{or}
\\
 \exists k \geq 1 \ \exists i_1, \ldots i_k \in I \ : \ a_{i, i_1}\;, 
a_{i_1, i_2}\;,
\;\ldots\;,
a_{i_k, j} < 0\;,
\end{array} \right. 
\end{align*}
and let $(I_\alpha)_{\alpha} \subseteq I$ denote the corresponding equivalence classes.
Then the following identities hold:
\begin{align}
\Ker A & = \{ (x_i)  \in \R^n \ | \ 
x_i = x_j \text{ whenever } i \sim j \}\;,\label{eq:Ker}\\
\label{eq:Ran}
\Ran A & = \Big\{ (x_i) \in \R^n \ | \ \forall \alpha \ : \
 						    \sum_{i \in I_\alpha}  x_i = 0	\Big\}\;.
\end{align}
\end{lemma}

\begin{proof}
First we remark that the assumptions $(1)$ -- $(3)$ imply that $a_{ij} = 0$ if $i \in I_\alpha$ and $j \in I_\beta$ for some $\alpha \neq \beta$.
Furthermore, 
it suffices to show \eqref{eq:Ker}, since \eqref{eq:Ran} then follows by duality. 

To show ``$\supseteq$'', suppose that $x = (x_i)$ satisfies $x_i = x_j$ whenever $i \sim j$. Fix $k \in I$ and take $\beta$ such that $k \in I_\beta$.
Using the remark and $(3)$, it follows that
\begin{align*}
\sum_{j \in I} a_{kj} x_j  
= \sum_{j \in I_\beta} a_{kj} x_j 
= x_k \sum_{j \in I_\beta} a_{kj}
= x_k \sum_{j \in I} a_{kj}
= 0\;,
\end{align*}
which yields the desired inclusion.

Conversely, to show ``$\subseteq$'', we use the identity 
\begin{align*}
 2 x_{ij} = x_i^2 + x_j^2 - (x_i - x_j)^2
\end{align*}
to write, for $x = (x_i)$,
\begin{align*}
2 \ip{A x, x}
& = 2\sum_{i,j \in I} a_{ij} x_i x_j
\\& =  \sum_{i \in I} x_i^2 \sum_{j \in I} a_{ij}
   +\sum_{j \in I} x_j^2 \sum_{i\in I} a_{ij} 
   -\sum_{i,j \in I} a_{ij} (x_i - x_j)^2\;.
\end{align*}
Using (3) and the symmetry of $A$ we infer that
\begin{align*}
 \ip{A x, x} = - \frac12\sum_{i,j\in I} a_{ij} (x_i - x_j)^2\;.
\end{align*}
Consequently, if $Ax = 0$, it follows that $\ip{Ax,x} = 0$, hence $x_i = x_j$ whenever $i \sim j$, which completes the proof.
\end{proof}

\section{Uniqueness of the metric on the two-point space}
\label{sec:proof}
In this appendix we shall prove Proposition \ref{prop:metric-uniqueness}.
First we need two definitions. Let $(M,d)$ be a metric space.

\begin{definition}\label{def:}
Let $I \subseteq \R$ be an interval and let $1 \leq p < \infty$. A curve $\gamma : I \to M$ is said to be \emph{$p$-absolutely continuous} if there exists a function $m \in L^p(I;\R)$ such that
\begin{align*}
  d( \gamma(s), \gamma(t)) \leq \int_s^t m(r) \dd r
\end{align*}
for all $s, t \in I$ with $s \leq t$. The curve $\gamma$ is \emph{locally $p$-absolutely continuous} if it is $p$-absolutely continuous on each compact subinterval of $I$.
\end{definition}

We shall use the notation
$
 \gamma \in AC^p(I;M)
$
and 
$ 
 \gamma \in AC_{\rm loc}^p(I;M)
$
respectively.

The following notion of gradient flow in a metric space $(M,d)$ has been studied in great detail in \cite{AGS08}.
  \begin{definition} \label{def:EVI}
Let $\cF : M \to \R \cup \{ + \infty\}$ be lower-semicontinuous and not identically $+ \infty$. A curve $\gamma \in C([0,\infty);M) \cap AC_{\rm loc}^2((0,\infty);M)$
is said to satisfy the  \emph{evolution variational inequality} $(\rm{EVI}_\lambda(\cF))$ if, for any $y \in
\Dom(\cF),$ the inequality 
 \begin{align} \label{eq:EVI}
 \frac{1}{2}\ddt d^2(\gamma(t), y) + \frac{\lambda}{2}
d^2(\gamma(t),y) \leq \cF(y) - \cF(\gamma(t)) 
\end{align}
holds a.e. on $(0,\infty)$.
 \end{definition}
   
\begin{proof}[Proof of Proposition \ref{prop:metric-uniqueness}]
Recall that $\bar\beta = \frac{p-q}{p+q}$.
Let $\beta \in (\bar\beta,1)$ and suppose that there exists $\alpha \in (-1,1)$ such that 
\begin{align} \label{eq:geod-eq}
  \cM(\rho^{\bar\beta}, \rho^\beta)  = \cM(\rho^{\bar\beta}, \rho^\alpha) + \cM(\rho^\alpha, \rho^\beta)\;.
\end{align}
We claim that $\alpha \in [\bar\beta,\beta]$. 
To prove this, suppose first -- to obtain a contradiction -- that $\alpha > \beta$. Then there exists $T > 0$ such that $e^{T (K-I)} \rho^\alpha = \rho^\beta$, hence \eqref{eq:EVI} implies that 
\begin{align*}
 \cM(\rho^\beta, \rho^{\bar\beta})^2 - \cM(\rho^\alpha, \rho^{\bar\beta})^2 
   \leq 2T (\cH(\rho^{\bar\beta}) - \cH(\rho^{\beta})) 
   \leq 0\;.
\end{align*}
In view of \eqref{eq:geod-eq}, it follows that $\cM(\rho^\alpha, \rho^\beta) = 0$, thus $\alpha = \beta$, which contradicts the assumption.
Suppose now that $\alpha < \bar\beta$. Adding \eqref{eq:geod-eq} and the inequality in \eqref{item:glueing} we infer that $\rho^\alpha = \rho^{\bar\beta}$, hence $\alpha = \bar\beta$, which proves the claim.

Now, fix $\beta \in (\bar\beta,1)$ and let $t \mapsto \rho^{\psi(t)}$ be a speed-1 geodesic with $\psi(0) = \bar\beta$ and $\psi(T) = \beta$ where $T = \cM(\rho^{\bar\beta},\rho^\beta)$. For $0 \leq s < t \leq T$ we then have $\cM(\rho^{\bar\beta},\rho^{\psi(t)}) = \cM(\rho^{\bar\beta},\rho^{\psi(s)}) +\cM(\rho^{\psi(s)},\rho^{\psi(t)})$, thus the claim implies that $\psi(s) \leq \psi(t)$. Since $\psi$ is a geodesic, we have $\psi(s) \neq \psi(t)$, thus $\psi$ is strictly increasing on $[0,1]$.

Now we claim that $\psi$ is continuous on $[0,T]$. To show this, take $t \in (0,T)$. Since $\psi$ is increasing, the limits $\psi(t-)$ and $\psi(t+)$ exist and for any $\eps > 0$ we have $\cM(\rho^{\psi(t-)}, \rho^{\psi(t+)} )\leq \cM(\rho^{\psi(t - \eps)}, \rho^{\psi(t - \eps)}) = 2\eps$, thus $\psi(t-) = \psi(t+)$. A similar argument shows that $\psi$ is continuous at $0$ and $T$, thus $\psi$ is continuous on $[0,T]$. 
Since $\psi$ is continuous and strictly increasing we infer that the mapping $\psi : [0,T] \to [\bar\beta,\beta]$ is surjective.
As a consequence, the inverse mapping $\phi : [\bar\beta,\beta] \to [0,T]$ is well-defined, and continuous and strictly increasing as well.

Note that the mapping 
\begin{align} \label{eq:isometry}
  I :  t \mapsto \rho^{\psi(t)}
\end{align}
defines an isometry from $[0,T]$ endowed with the euclidean metric onto $\{ \rho^\alpha : \alpha \in [0,\beta]\}\subseteq \cP_*(\cX)$ endowed with the metric $\cM$. The inverse mapping is given by
\begin{align*}
 J : \rho^\alpha \to \phi(\alpha)\;.
\end{align*}
Since $u : t \mapsto \rho^{\beta_t}$ is a $2$-absolutely continuous curve satisfying $\text{EVI}_0(\cH)$ for the metric $\cM$, \eqref{eq:isometry} implies that the mapping 
\begin{align*}
t \mapsto \tilde u(t) := J(u(t)) = \phi(\beta_t)
\end{align*}
is a $2$-absolutely continuous curve satisfying $\text{EVI}_0(\tilde \cH)$ where $\tilde\cH := \cH \circ I$, for the euclidean metric.
It follows that the mapping $\phi : [\bar\beta,\beta] \to [0,T]$ itself is absolutely continuous, hence almost everywhere differentiable, and the same holds for its inverse $\psi$. Moreover, the identity
\begin{align} \label{eq:chain-psi-phi}
 \psi'(\phi(\alpha)) \phi'(\alpha) = 1
\end{align}
holds for a.e. $\alpha \in [\bar\beta,\beta]$.

For any $\alpha \in [\bar\beta, \beta]$ we have
\begin{align*}
\tilde \cH(\phi(\alpha)) 
  & = \cH(I(\phi^\alpha))
  = \cH(\rho^{\alpha})
  \\& = \frac{q}{p+q}
    f\Big( \frac{p + q}{q} \frac{1- \alpha}{2}   \Big) 
    + \frac{p}{p+q}
    f \Big(  \frac{p + q}{p} \frac{1+ \alpha}{2} \Big)\;,
\end{align*}
thus, for $r \in (0,T)$,
\begin{align*}
  \tilde\cH(r) =  \frac{q}{p+q}
    f\Big( \frac{p + q}{q} \frac{1- \psi(r)}{2}   \Big) 
    + \frac{p}{p+q}
    f \Big(  \frac{p + q}{p} \frac{1+ \psi(r)}{2} \Big)\;.
\end{align*}
It follows that $\tilde\cH$ is a.e. differentiable and 
the identity
\begin{align}\label{eq:Htilde-prime}
 \tilde\cH'(r) = \frac{\psi'(r)}2 
   \Big[ f'\Big(  \frac{p + q}{p} \frac{1+ \psi(r)}{2} \Big)
    -  f'\Big( \frac{p + q}{q} \frac{1- \psi(r)}{2}   \Big)  \Big]
\end{align}
holds a.e.

Since $t \mapsto \tilde u(t)$ is a $2$-absolutely continuous curve satisfying $\text{EVI}_0(\tilde\cH)$ and since the functional $\tilde\cH$ is differentiable a.e., it follows from \cite[Proposition 1.4.1]{AGS08} that the gradient flow equation
\begin{align*}
  \tilde u'(t) = -  \tilde\cH'(\tilde u(t))
\end{align*}
holds almost everywhere.

Since $\phi$ is differentiable a.e., the left-hand side equals a.e.
\begin{align*}
  \tilde u'(t) 
  = \ddt \phi(\beta_t) 
  =   \big(p(1- \beta_t) -  q(1 + \beta_t)\big)\phi'(\beta_t)\;.
\end{align*}
Taking \eqref{eq:chain-psi-phi} into account, it follows from \eqref{eq:Htilde-prime} that the right-hand side equals a.e.
\begin{align*}
 \tilde\cH'(\tilde u(t)) = \frac{1}{2 \phi'(\beta_t)} 
 \Big[  f'\big( \rho^{\beta_t}(b) \big)
 				-  f'\big( \rho^{\beta_t}(a) \big) \Big]
\end{align*}
Combining the latter two inequalities we infer that for a.e. $\alpha \in [\bar\beta,\beta]$,
\begin{align*}
 \big(q(1 + \alpha) - p(1- \alpha) \big)\phi'(\alpha)
     = \frac{1}{2 \phi'(\alpha)} 
 \Big[ f'\big( \rho^{\alpha}(b) \big)
 				-  f'\big( \rho^{\alpha}(a) \big) \Big]
\end{align*}
Since $\phi$ is absolutely continuous,
\begin{align*}
 \phi(\beta)
   = \int_{\bar\beta}^\beta \phi'(\alpha) \dd \alpha
   = \int_{\bar\beta}^\beta \sqrt{\frac{f'\big( \rho^{\alpha}(b) \big)
 				-  f'\big( \rho^{\alpha}(a) \big) }
				{2 \big(q(1 + \alpha) - p(1- \alpha)\big)}}
     \dd \alpha\;.
\end{align*}
hence, since $t \mapsto \psi(t)$ is a geodesic, we obtain for $\bar\beta < \alpha < \beta$,
\begin{align*}
 \cM(\rho^\alpha, \rho^\beta)
   =  \cM(\rho^{\psi(\phi(\alpha))}, \rho^{\psi(\phi(\beta))})
   = C (\phi(\beta) - \phi(\alpha) )\;.
\end{align*}
Thus the distance between $\rho^\alpha$ and $\rho^\beta$ is uniquely determined for all $\alpha, \beta \geq \bar\beta$. The same argument shows that the distance is uniquely determined for $\alpha, \beta \leq \bar\beta$. The case $\alpha < \bar\beta < \beta$ follows from the assumption \eqref{item:glueing}.
\end{proof}

\bibliographystyle{ams-pln}
\bibliography{cube}
 \end{document}